\documentclass[12pt]{article}

\usepackage[latin1]{inputenc} 
\usepackage{latexsym}
\usepackage{graphicx}         
\usepackage{epsfig}           
\usepackage{amsmath}          
\usepackage{amsfonts}
\usepackage{mathrsfs}
\usepackage{fancyhdr}

\DeclareGraphicsExtensions{.eps,.ps,.pdf,.bmp}

\newcommand{\A}{\mathrm{A}}

\newcommand{\K}{\mathrm{K}}

\newcommand{\La}{\Lambda}

\newcommand{\Om}{\Omega}

\newcommand{\Ga}{\Gamma}

\newcommand{\sig}{\sigma}

\def \rat{ {\rm Q}\kern-.65em {}^{{}_/ }}

\newcommand {\p} {\partial}

\newcommand {\hs} {\hspace}

\newtheorem {Def} {\noindent Definition}[section]

\newtheorem {thm} {\noindent Theorem}[section]
\newtheorem {pr} {\noindent Proposition}[section]

\newtheorem {Le} {\noindent Lemma}[section]
\newtheorem {co} {\noindent Corollary}[section]

\makeatletter
\renewcommand{\section}{\@startsection {section}{1}{0pt}
    {-3.5ex \@plus -1ex \@minus -.2ex}
    {2.3ex \@plus.2ex}
    {\normalfont\normalsize\bfseries}}
\newcommand{\SUMMARY}{\@startsection {section}{1}{0pt}
    {-3.5ex \@plus -1ex \@minus -.2ex}
    {2.3ex \@plus.2ex}
    {\normalfont\normalsize\centering\bfseries}*{{\bf R\'esum\'e}}}
    \makeatother

\begin{document}
\vspace*{2cm} \normalsize \centerline{\Large \bf Riesz bases for $L^2(\p \Om)$ and regularity }
\vspace*{0.5cm} \normalsize \centerline{\Large \bf  for the Laplace equation in Lipschitz domains }
\vspace*{1cm}

\centerline{\bf {\underline{Abdellatif CHA\"{I}RA} } \footnote{ Corresponding author: a.chaira@fs.umi.ac.ma}, \underline{Soumia TOUHAMI } \footnote{s.touhami@edu.umi.ac.ma}}

\vspace*{0.5cm}


\vspace*{2cm}
 Universit\'{e} Moulay Ismail,  Facult\'{e} des Sciences, \\
  Laboratoire de Math\'{e}matiques et leurs Applications, \'{E}quipe EDP et Calcul Scientifique,\\
  BP $11201$ Zitoune, 50070 Mekn\`{e}s, Maroc. \\ \\
\vspace*{1cm}

\noindent {\bf Abstract.}
In a paper from 1996,  D. Jerison and C. Kenig among other results provided  a $H^{1/2}$ regularity result for the Dirichlet problem for the Laplace equation in Lipschitz domains. In this article, we adopt a Hilbertian approach to construct two Riesz bases for $L^2(\p \Om),$ which will allow to find in a different way some of the results of D. Jerison and C.  Kenig, and G. Savar\'{e} (1998)  about the regularity issue of the Laplace equation. 

\vspace*{0.5cm}


\vspace*{0.5cm}

\setcounter{equation}{0}
\section{Introduction and results}

Given a bounded Lipschitz domain $\Omega \subset  \mathbb R^d, d\geq 2,$ with boundary $\p \Omega$, let us consider for $g$ defined on $\p \Om,$ the  Dirichlet problem for the Laplace equation 
 \begin{equation}
            \begin{cases}
              \Delta v=0 & \text{  }  (\Omega) \\
             v=g & \text{}  (\partial \Omega).
            \end{cases}
   \end{equation}
   
We briefly recall the history of the problem (1). If $g$ is continuous on $\p \Om$ it is well known that $\Omega$ being regular for the Laplacian  $\Delta,$ the problem (1) has a unique solution given by
   $$v(x) = \int_{\p \Omega} g(y) \ d \omega ^x(y),$$
   where $\omega ^x $ is the harmonic measure for $\Omega$ with  pole $x \in \Omega.$  For $g\in H^{1/2}(\p \Omega),$  the problem $(1)$ is variational and has a unique solution according to the Hilbertian theory of Sobolev spaces \cite{Ga}. \\
   
   \ \ \ \ \ \ The problem (1) when the data consisted either of functions in $L^{2}(\p \Omega)$ or of functions with first derivatives in $L^2(\p \Omega),$ had  attracted significant  research attention. This began with the work of  J. Ne\v{c}as \cite{N}. Using {\bf  Rellich Identity}: 
    \begin{equation*} \int_\Omega (m,\nabla u)_{\mathbb{R}^d}\Delta udx=
   -\sum_{i,j=1}^d\int_\Omega
   \frac{\partial m_i}{\partial x_j}\frac{\partial u}{\partial x_i}\frac{\partial u}{\partial x_j}dx
   + \ \frac{1}{2}\int_\Omega \mbox{div}(m)|\nabla u|^2dx \  + \frac{1}{2}\int_{\partial
   \Omega}(m,\nu)_{\mathbb{R}^d}|\partial_\nu u|^2d\sigma,
   \end{equation*}
   where  $m\in ({\cal C}^\infty(\mathbb R^d))^d$ is a vector field, $\partial_\nu$ is the normal derivative operator associated to $\Omega$ and $(.,.)_{\mathbb{R}^d}$ denotes the inner product on $\mathbb R^d,$ J. Ne\v{c}as proved the following result
     \begin{itemize}
     \item  {\bf  Rellich-Ne\v{c}as lemma.} {\it Let $\Omega \subset \mathbb{R}^d$  be a bounded Lipschitz domain. Then, there exists a constant  $c_{\Omega}>0$ depending on the geometry of $\Omega$ such that for all $u\in H^1_\Delta(\Omega)\cap H^1_0(\Omega)$ }
   $$
   \|\partial_\nu u\|_{L^2(\partial\Omega)}\leq c_{\Omega} \ \|\Delta
   u\|_{L^2(\Omega)},
   $$
   where $H^1_\Delta(\Omega) =\{ \ \  u\in H^1(\Omega) \ \  | \ \  \Delta u \in L^2(\Omega) \ \  \}$.
     \end{itemize}
   The Rellich-Ne\v{c}as lemma allows to define the very weak solution of the Dirichlet problem for the Laplace equation (1). Indeed, we say that $v\in L^2(\Om)$ is a \textbf{very weak solution} of the problem (1) if for all $u\in H_{\Delta}^1(\Om) \cap H_0^1(\Om)$ we have 
   \begin{equation}
   - \int_{\Om} v \Delta u \ dx + \int_{\p \Om} g \  \partial _{\nu} u \ d\sigma  =0 .
   \end{equation}
   
   The last formulation makes sens according to Rellich-Ne\v{c}as lemma and means that 
   $$ A'v= \mu,$$ where $A'$ from $ L^2(\Om)$ to the dual space of $H_{\Delta}^1(\Om) \cap H_0^1(\Om)$ denoted $   (H_{\Delta}^1(\Om) \cap H_0^1(\Om))^{\prime},$  is the transpose of the Laplacian with Dirichlet conditions
   $$\varphi \longrightarrow A\varphi= -\Delta \varphi$$ defined 
   from $ H_{\Delta}^1(\Om) \cap H_0^1(\Om)$ to $L^2(\Om),$ and where $\mu$ is the linear form given by
   $$\varphi \longrightarrow 
    -\int_{\p \Om} g \  \partial_{\nu} \varphi \ d\sigma ,$$
   which is continuous according to Rellich-Ne\v{c}as lemma. Since $A$ is an isomorphism from $ H_{\Delta}^1(\Om) \cap H_0^1(\Om)$ into $L^2(\Om),$ it follows that $A^{\prime}$ is also an isomorphism from $L^2(\Om)$ to $(H_{\Delta}^1(\Om) \cap H_0^1(\Om))^{\prime},$ and this proves the existence and the uniqueness of $v\in L^2(\Om),$  solution of (1). Dahlberg in \cite{Da}, established that the harmonic measure  and the surface measure associated to $\Omega$  are mutually absolutely continuous, furthermore, the Random-Nikodym derivative of harmonic measure with respect to surface measure satisfies a reverse H\"{o}lder inequality which allows to solve the problem (1) with data in $L^2(\p \Om).$  In \cite{JK1}, D. Jerison and C. Kenig provided another proof of Dahlberg's results using an integral identity due to Rellich, and after that in \cite{JK2},  they gave optimal estimates for the Dirichlet problem when the data has one derivative in $L^2(\p \Om),$ where they combined Rellich formulas with Dahlberg's results. D. Verchota in \cite{V}, following the works of Coifman-McIntosch and Meyer \cite{CMM}, had been interested to the invertibility of classical layer potentials for Laplace equation on the boundaries of bounded Lipschitz domains and the applications to the Dirichlet and Neumann problems. In \cite{JK},  D. Jerison and C. Kenig   studied the inhomogenous Dirichlet problem for the Laplacian in Lipschitz domains with data in trace spaces, where they used the strategy of reduction to the homogenous problem. The two important tools in their paper were the investigation of traces of Sobolev spaces on the boundary and the characterization of Sobolev and Besov spaces of harmonic functions.  Savar\'{e} in \cite{S}, developped a variational argument based on the usual Niremberg's difference quotient technique to deal with the regularity of the solutions of Dirichlet and Neumann problems for some linear and quasilinear elliptic equations in Lipschitz domains.\\
   
\ \ \ \ \ \ \ The main purpose of this paper is to construct two Riesz bases for $L^2(\p \Om)$ (see \S 4) and show how it will be possible to give another proof of  the $H^{1/2}$ regularity results about the Dirichlet problem for the Laplacian previousely established by Jerison and Kenig in \cite{JK2} and by Savar\'{e} in \cite{S} (see \S 5). In the following we give a first description of the approach that we will follow in this paper and which will be detailed in the next sections.\\

 Consider the solution operator of the problem (1)  $$ K: L^2(\p \Omega) \longrightarrow L^2(\Omega) ,\ \   g \mapsto v,$$
where $v$ is the very weak solution of (1) and consider its adjoint operator $K^*,$ which takes each $f \in L^2(\p \Omega)$ to $-\partial_{\nu} u^0$ into $L^2(\p \Omega),$ where $u^0$ is the solution of the Dirichlet problem for the following  Poisson equation 
\begin{equation}
              \begin{cases}
                -\Delta u^0=f & \text{  }  (\Omega) \\
                \Gamma u^0=0 & \text{}  (\partial \Omega),
              \end{cases}
              \end{equation}

where $\Gamma$ is the trace operator from $H^1(\Omega)$ to $L^2(\partial \Omega).$ Consider also the embedding operator from $H^1(\Omega)$ into $L^2(\Omega)$ denoted $E.$ For $f\in L^2(\Omega),$ the adjoint operator $ E^*$ is the solution operator of Robin problem for the following Poisson equation
  \begin{equation}
             \begin{cases}
               -\Delta u=f & \text{  }  (\Omega) \\
              \partial_{\nu} u + \Gamma u=0 & \text{}  (\partial \Omega).
             \end{cases}
             \end{equation}

By setting $E_1^*= E^*-E_0^*$ and $u^1 =E_1^*f,$ where $E_0^*$ denotes the solution operator of (3), it follows that $u^1$ is a solution of the following Dirichlet problem for the Laplace equation
 \begin{equation}
              \begin{cases}
                -\Delta u^1=0 & \text{  }  (\Omega) \\
               \Gamma u^1=\Gamma u & \text{}  (\partial \Omega),
              \end{cases}
              \end{equation}
   where  $u$ is the solution of (4).  
   Let us set $\Gamma_0^*= F_1^*(I+F_1F_1^*)^{-1/2} \Gamma^*,$ where $F_1$ is the Moore-Penrose inverse of the adjoint operator $E_1 =(E_1^*)^*,$ and $\Gamma^*$ is the adjoint of the trace operator $\Gamma.$ We will show in section 3 of this paper that
   $$ \Gamma_0^* K^* = (I+F_1^*F_1)^{-1/2} P_{\mathcal H(\Omega)}  $$ is compact and self-adjoint, where $P_{\mathcal H(\Om)}$ is the orthogonal projection onto the space of harmonic square-integrable functions which is called the Bergman space and denoted in this text by $\mathcal H(\Omega).$ Consequently, there exists a sequence of couples $((\kappa_n,\phi_n))_{n\geq1}  \in \mathbb R _+^* \times \mathcal H(\Omega) $ associated with $\Gamma_0^* K^*$ such that $$\Gamma_0^* K^* \phi_n =\kappa_n^2 \phi_n.$$
Moreover, $(\phi_n)_{n\geq 1}$ is an orthonormal basis for $\mathcal H(\Om).$ By setting for all $n\geq 1,$
   $$ \Gamma_0^* \phi_n = \kappa_n y_n \ \mbox{and} \ K^* \phi_n = \kappa_n g_n,$$ 
   the main purpose of the present work is to prove the following result.
\begin{thm}
   Let $\Omega \subset \mathbb R^d$ be a bounded Lipschitz domain. Then, the sequences $(g_n)_{n\geq 1}$ and $(y_n)_{n\geq 1}$ defined above, are Riesz bases for $L^2(\p \Om).$
   \end{thm}
    
 One of the main consequences of Theorem 1.1 is the following classical regularity result. 

   \begin{thm} Let $\Omega \subset \mathbb R ^d$ be a bounded Lipschitz domain. Then, for $g\in L^2(\partial\Omega)$, the very weak solution of the Dirichlet problem for the Laplace equation (1) lies in $H^{1/2}(\Omega)$ and there exist two positive constants $c_{\Omega}$  and $c_{\Omega}^{\prime}$ depending on the geometry of $\Omega$ such that 
   $$
   c_{\Omega}^{\prime} \ \|g\|_{L^2(\partial\Omega)} \ \leq \|v\|_{H^{1/2}(\Omega)}\leq c _{\Omega}\ \|g\|_{L^2(\partial\Omega)},
   $$ 
   Moreover, the solution operator $K$ is compact and injective.
   \end{thm}

   The plan of the paper is the following: the next section contains some known and new facts about the Moore-Penrose inverse and a brief recall of some preliminary results for Riesz bases and related sequences, and also some basic results for  Sobolev spaces in Lipschitz domains. In  section 3, we present the main key tools to deal with Theorem 1.1 and Theorem 1.2. Section 4 will be devoted to study  the sequences $(g_n)_{n\geq 1}$ and $(y_n)_{n\geq 1}$ and to present the remaining arguments to conclude our main result (Theorem 1.1). In section 5, a regularity result for the Dirichlet problem for the Laplace equation (1) will be derived (Theorem 1.2).

   \section{Preliminaries and basic results}
   
 Let $(\mathcal H_1 , (.,.)_1)$ and  $(\mathcal H_2 , (.,.)_2)$  be two Hilbert spaces  with the associated inner products $ (.,.)_1, (.,.)_2 $ and the induced norms  $\|.\|_1, \|.\|_2,$  and throughout this article, unless otherwise mentioned, they will be simply denoted $\mathcal H_1$ and $\mathcal H_2.$ A linear operator from $\mathcal H_1$ to $\mathcal H_2$ is a pair consisting of a subspace $\mathcal D(A)$ of $\mathcal H_1$ together with a linear map $A: \mathcal D(A) \longrightarrow \mathcal H_2.$ We call $\mathcal D(A)$ the domain of the operator $A$ and write  $(A,\mathcal D(A))=A.$   $\mathcal N(A)$ denotes its null space, $\mathcal R(A)$ its range space and $\mathcal G(A)$ its graph. In the case $(A,\mathcal D(A))$ is bounded, we write simply $A.$   The set of all bounded operators from $\mathcal H_1$ into $\mathcal H_2$  is denoted by $\mathcal B(\mathcal H_1, \mathcal H_2),$ and if $\mathcal H_1 = \mathcal H_2,$  $\mathcal B(\mathcal H_1, \mathcal H_2)$ is denoted $ \mathcal B(\mathcal H).$
  For two linear operators $(A,\mathcal D(A))$ and $(B,\mathcal D(B))$  from $\mathcal H_1$ into $\mathcal H_2,$  $(B,\mathcal D(B))$  is called an extension of $(A,\mathcal D(A))$ if
        $$ 
        \begin{array}{ccccc}
       \mathcal  D(A)\subset  \mathcal D(B) & \mbox{and} \\
        \forall x\in \mathcal D(A), Ax  =  Bx,
        \end{array}
        $$
   and this fact is denoted by  $A\subset B.$ \\ For a linear operator $(A,\mathcal D(A))$  on a Hilbert space $\mathcal H,$ there are several ways of defining the notion of positivity, in this paper, this corresponds to the following 
     $$ (Ax,x) \geq 0  \ \ \forall  x \in \mathcal D(A) ,$$
      in such case we write $A\geq 0$ and say that $(A,\mathcal D(A))$ is positive.  $(A,\mathcal D(A))$ is said to be densely defined if $\mathcal D(A)$ is dense in $\mathcal H_1, i.e., \overline{\mathcal D(A)} = \mathcal H_1,$ where $\overline{\mathcal D(A)}$ denotes the closure of $\mathcal D(A).$  $(A,\mathcal D(A))$ is said to be closed if its graph is closed in $\mathcal H_1 \times \mathcal H_2,$   where the inner product in  $\mathcal H_1 \times \mathcal H_2$ is defined for all $x,u \in \mathcal H_1$ and $y,v \in \mathcal H_2$ by $$((x,y),(u,v))= (x,u)_1 + (y,v)_2.$$

     The set of all closed densely defined operators from $\mathcal H_1$ into $\mathcal H_2$ is denoted by $\mathcal C(\mathcal H_1, \mathcal H_2).$

The adjoint of a densely defined operator from $\mathcal H_1$ into $\mathcal H_2$ is denoted $(A^*,\mathcal D(A^*))$ where
     $\mathcal D(A^*)$ is defined to be the set of all $y \in \mathcal H_2$ for which there exists $z\in \mathcal H_1$ such that
          $$(Ax,y)_2 =(x,z)_1 \ \
          \forall x\in \mathcal D{(A)}.$$ Since $\mathcal D(A)$ is dense, it follows that $z$ is unique. We put $A^*y =z,$ then we have: 
   $$(x,A^*y)_1= (Ax,y)_2, \forall x \in \mathcal D( A), y\in \mathcal D(A^*),$$
  and $A^*$ is closed. Moreover, if $A$ is closed,  $A^*$ is densely defined.
   \begin{Le}
   Let $\mathcal H_1$  and $\mathcal H_2$ be two Hilbert spaces and  $(A, \mathcal D(A)), (B, \mathcal D(B))$ be two linear operators from $\mathcal H_1$ into $\mathcal H_2$ such that $A\subset B.$ Then if $\mathcal D(A)$ is dense, we have $B^* \subset A^*.$ 
   \end{Le}
     A linear operator $ (A, \mathcal D(A))$ on a Hilbert space $\mathcal H$ is said to be self-adjoint if $A^*=A$ which means that $\mathcal D(A^*) = \mathcal D(A)$ and that $A^*x= A x$ for all $x \in \mathcal D(A).$ 
       Many of the operators which we shall study in this paper are positive self-adjoint and the condition of self-adjointness is of profound importance to define the powers of any fractional order of $(A,\mathcal D(A)).$ 
    A bounded linear operator $A$ from $\mathcal H_1 $ to $\mathcal H_2$ is said to be compact if for any bounded sequence $(f_n)_{n\geq1} $  of elements of  $ \mathcal H_1,$ the sequence $(A f_n)_{n\geq1}$ has a norm convergent subsequence. The following theorem is stated in \cite[Theorem 3.4]{Co1}.
   \begin{thm} (\textbf{Schauder's Theorem}) Let $\mathcal H_1$  and $\mathcal H_2$ be two Hilbert spaces and $A \in \mathcal B(\mathcal H_1,\mathcal H_2).$ Then, $A$ is compact if and only if its adjoint $A^*$ is compact.
      \end{thm}
  For further lectures, see \cite{Co} and \cite{K}.\\

    When an operator is not invertible in the strict sense, one can define its Moore-Penrose inverse. The next subsection is devoted to provide some known and new  facts about this concept that will play a key role in this text.

    \subsection{The Moore-Penrose Inverse}

Let  $ {\cal H}_1, {\cal H}_2$ be two Hilbert spaces, $A \in {\cal C}({\cal H}_1, {\cal H}_2) $ a closed densely defined operator and  $(A^*,\mathcal D(A^*))$  its adjoint. The Moore-Penrose inverse  of   $(A,\mathcal D(A)) $ denoted  $(A^\dagger,{\cal D}(A^\dagger)) $ is defined as the unique linear operator in  ${\cal C }({\cal H}_2, {\cal H}_1)$  such that
       $$
       {\cal D}(A^\dagger)={\cal R}(A)\oplus {\cal N}(A^*), \  \mathcal N(A^{\dagger })= \mathcal N(A^*) $$ 
      and satisfying the followings
       
        \begin{equation*}
              \begin{cases}
                AA^\dagger A=A & \text{  }  \\
                A^\dagger AA^\dagger=A^\dagger & \text{}
              \end{cases}
             \hs{1cm}   \hs{1cm}
                    \begin{cases}
                      AA^\dagger \subset P_{\overline{{\cal R}(A)}} & \text{  }  \\
                     A^\dagger A \subset P_{\overline{{\cal R}( A^\dagger)}} & \text{},
                    \end{cases}
              \end{equation*}
where $ P_{\overline{{\cal R}(A)}}$ and $ P_{\overline{{\cal R}( A^\dagger)}}$ denote the orthogonal projections onto $\overline{{\cal R}(A)}$ and $\overline{{\cal R}( A^\dagger)}$ respectively.
       Moreover, $(A, \mathcal D(A))$ is the Moore-Penrose inverse of $(A^\dagger,{\cal D}(A^\dagger))$ and ${\cal R}(A)$ is closed  if and only if $(A^\dagger,{\cal D}(A^\dagger))$ is bounded. According to a fundamental result of Von Neumann (see \cite{Gro} and \cite{L}), for $A\in \mathcal C(\mathcal H_1,\mathcal H_2),$ the operators $(I+AA^*)^{-1}$ and $A^*(I+AA^*)^{-1}$ are everywhere defined and bounded. Moreover, $(I+AA^*)^{-1}$ is self-adjoint.
      Similarly, the operators $(I+A^*A)^{-1}$ and $A(I+A^*A)^{-1} $ are everywhere defined and bounded, and  $(I+A^*A)^{-1}$ is self-adjoint. Moreover, we have the following
    
      $$(I+AA^*)^{-1} A \subset A(I+A^*A)^{-1}$$
      and
      $$(I+A^*A)^{-1} A^*  \subset A^*(I+AA^*)^{-1}.$$
    
      (see \cite{Gro} and \cite{L}).\\ \\
    In the following, we state some identities that go back to Labrousse \cite {L}:
      \begin{pr}  Let $A \in \mathcal C(\mathcal H_1, \mathcal H_2)$ and $B \in \mathcal C(\mathcal H_2,\mathcal H_1)$ such that $ B=A^{\dagger}$, then
      \begin{enumerate}
      \item  $A(I+A^*A)^{-1} = B^*(I+BB^*)^{-1};$
      \item $(I+A^*A)^{-1}+(I+BB^*)^{-1}= I+P_{\mathcal N(B^*)};$

       \item  $A^*(I+AA^*)^{-1} = B(I+B^*B)^{-1};$
      \item $(I+AA^*)^{-1}+(I+B^*B)^{-1}= I+P_{\mathcal N(A^*)};$
       \item $(I+AA^*)^{-1}+(I+B^*B)^{-1}= I$  (if $A^*$ is injective);
    \item $\mathcal N(A^*(I+AA^*)^{-1/2} ) =\mathcal N(A^*) = \mathcal N(B).$
      \end{enumerate}
     
     \end{pr}
   
    Some of the results  we will present in the rest of this subsection about the Moore-Penrose inverse, are stated for the first time and will prove useful throughout the rest of this paper.
     \begin{pr} Let  $ {\cal H}_1, {\cal H}_2$ be two Hilbert spaces, $A \in {\cal B}({\cal H}_1, {\cal H}_2) $ and $B$ its Moore-Penrose inverse, then
     for all $x\in \mathcal H_1$ one has 
     $$ \|x\|_1^2= \|B^*(I+BB^*)^{-1/2} x\|_2^2 + \|(I+BB^*)^{-1/2} x\|_1^2.$$
      Moreover, if $x\in \mathcal R(B)$ then
     $$\|x\|_1^2 = \|(I+BB^*)^{-1/2} x\|_1^2 + \|(I+A^*A)^{-1/2} x\|_1^2.$$
     \end{pr}
     \begin{pre}
     The first part of the proposition was proved in \cite{LM}. Now, 
     for $x\in \mathcal R(B)=\mathcal N(B^*)^{\perp}$ where $\mathcal N(B^*)^{\perp}$ denotes the orthogonal complement of $\mathcal N(B^*),$ we have according to the fourth item of Proposition 2.1 that
     $$ (I+A^*A)^{-1}x+ (I+BB^*)^{-1}x = x,$$
     which implies that
     \begin{eqnarray*}
     \|x\|_1^2 &=& (x,x)_1=\big((I+A^*A)^{-1}x+(I+BB^*)^{-1}x,x\big)_1\\
     &=& \|(I+A^*A)^{-1/2}x\|_1^2 + \|(I+BB^*)^{-1/2}x \|_1^2. \  \square
     \end{eqnarray*}   
     \end{pre}  \\
We will extensively make use of the following  result:
    \begin{pr}   Let  $ {\cal H}_1, {\cal H}_2$ be two Hilbert spaces, $A \in {\cal B}({\cal H}_1, {\cal H}_2) $ and $B$ its Moore-Penrose inverse, then the operator $B^*(I+BB^*)^{-1/2}$ is bounded, has a closed range and its Moore-Penrose inverse is given by 
    $$
    T_B=B(I+B^*B)^{-1/2}+A^*(I+B^*B)^{-1/2}.
    $$
    Moreover, the adjoint operator of $T_B$ is  $T_{B^*},$ where 
    $$T_{B^*}=B^*(I+BB^*)^{-1/2}+A(I+BB^*)^{-1/2}.$$
    \end{pr}
    \begin{pre}
    For $x\in \mathcal H_1,$ we have according to Proposition 2.2 that
    $$ \|x\|_1^2= \|B^*(I+BB^*)^{-1/2} x\|_2 ^2+ \|(I+BB^*)^{-1/2} x\|_1^2,$$
     and if $x\in \mathcal R(B),$ then 
    $$\|x\|_1^2 = \|(I+BB^*)^{-1/2} x\|_1^2 + \|(I+A^*A)^{-1/2} x\|_1^2,$$
    which implies that for all $x\in \mathcal R(B)=\mathcal N(B^*)^{\perp},$ we have
    $$\|B^*(I+BB^*)^{-1/2} x\|_2 = \|(I+A^*A)^{-1/2} x\|_1.$$
    Since $A$ is bounded, it follows that $(I+A^*A)^{-1/2}$ is bounded, invertible and has a bounded inverse. Moreover, there exists a positive constant $c$ such that for all $x\in \mathcal H_1$
    $$c \ \|x\|_1 \leq \|(I+A^*A)^{-1/2}x \|_1 \leq \|x\|_1$$ 
    and if $x\in \mathcal R(B),$ 
    $$c \ \|x\|_1 \leq \|B^*(I+BB^*)^{-1/2}x \|_2 \leq \|x\|_1.$$ 
   We therefore deduce that $B^*(I+BB^*)^{-1/2}$ has a bounded Moore-Penrose inverse, and a direct verification leads to 
    
    $$T_{B} B^*(I+BB^*)^{-1/2} T_{B} = T_{B}$$
    and that $$ \ B^*(I+BB^*)^{-1/2} T_{B}B^*(I+BB^*)^{-1/2} = B^*(I+BB^*)^{-1/2}.$$
     Moreover, we have 
    $$T_B B^*(I+BB^*)^{-1/2}=B(I+B^*B)^{-1/2} T_{B^*} = P_{\mathcal R(B)},$$
and
   $$T_{B^*} B(I+B^*B)^{-1/2}=B^*(I+BB^*)^{-1/2} T_{B} = P_{\mathcal R(B^*)}.$$
    Therefore, $T_B$ is the Moore-Penrose inverse of  $B^*(I+BB^*)^{-1/2}.$  On the other hand, since $$(B(I+B^*B)^{-1/2})^* =B^*(I+BB^*)^{-1/2}$$ and  $$(A^*(I+B^*B)^{-1/2})^*=A(I+BB^*)^{-1/2},$$ we obtain that
    \begin {eqnarray*}
    (T_{B})^*&=&(B(I+B^*B)^{-1/2})^*+(A^*(I+B^*B)^{-1/2})^*\\
    &=& B^*(I+BB^*)^{-1/2}+A(I+BB^*)^{-1/2}\\
    &=&T_{B^*}.
    \end{eqnarray*}
    Hence, $$(T_B)^*=T_{B^*}. \  \square$$
    \end{pre}
    \begin{co}
    The operator $B^*(I+BB^*)^{-1/2}$ is an isomorphism from $\mathcal N(B^*)^{\perp}$ to $\mathcal R(B^*).$ 
     \end{co}
      \begin{co}
         The operator $T_B$ is an isomorphism from $\mathcal R(B^*)$ to $\mathcal N(B^*)^{\perp}.$  
          \end{co}
        
      The next result provides a decomposition for an arbitrary  bounded operator in terms of its Moore-Penrose inverse.
      \begin{pr}
      Let $\mathcal H_1$ and $\mathcal H_2$ be two Hilbert spaces, $A\in \mathcal B(\mathcal H_1,\mathcal H_2)$ and $(B, \mathcal D(B))$ its Moore-Penrose inverse. Then, we have the following decomposition
      $$A= (I+B^*B)^{-1/2} T_{B^*},$$
      where $T_{B^*}=B^*(I+BB^*)^{-1/2}+A(I+BB^*)^{-1/2}.$
      \end{pr}
      \begin{pre}
      A direct verification leads to
      
      $$(I+B^*B)^{-1/2}T_{B^*}=(I+B^*B)^{-1/2} \Big (B^*(I+BB^*)^{-1/2}+A(I+BB^*)^{-1/2}\Big ).$$
Moreover, since $$(I+B^*B)^{-1/2} B^* \subset B^* (I+BB^*)^{-1/2},$$
and that $$ B^*(I+BB^*)^{-1}=A(I+A^*A)^{-1}$$ from the third item of Proposition 2.1, 
 it follows that

\begin{eqnarray*}
(I+B^*B)^{-1/2}T_{B^*} &=&  B^*(I+BB^*)^{-1}+A(I+BB^*)^{-1} \\
   &=& A(I+A^*A)^{-1}+A(I+BB^*)^{-1}\\
      &=& A \Big( (I+A^*A)^{-1}+(I+BB^*)^{-1} \Big ).
      \end{eqnarray*}
      Moreover, for $x\in \mathcal N(A),$ we have:
      \begin {eqnarray*}
      (I+B^*B)^{-1/2}T_{B^*}x&=&  A \Big( (I+A^*A)^{-1}+(I+BB^*)^{-1} \Big ) x \\
      &=& A(2x)\\
      &=& 2Ax \\
      &=&0.
      \end{eqnarray*}
      For $x\in \mathcal{R}(B),$ it follows  according to Proposition 2.1 that $$(I+B^*B)^{-1/2}T_{B^*}x=Ax.$$
      Hence, for all $x\in \mathcal H_1= \mathcal N(A) \oplus \mathcal R(B),$ 
      $$ Ax= (I+B^*B)^{-1/2})T_{B^*}x.\ \ \square $$
      \end{pre}
      Further detailed results concerning the Moore-Penrose inverse concept could be found in (\cite{Gro},\cite{L} and \cite{LM}).
   Another important theoretical background in Functional Analysis that will be useful in this paper is Riesz bases concept and related sequences, and most of the basic  results that we will remind here are stated  in (\cite{Ch}, \cite{H} and \cite{Y}).

\subsection{Riesz bases and related sequences} 

 A sequence $(x_k)_{k\geq1}$ in a Hilbert space $\mathcal H$ is said to be complete if $$ \overline{ span}(x_k)_{k\geq 1} = \mathcal H, $$ and minimal if each element of the sequence lies outside the closed linear span of the others, i.e., $$x_j \not \in \overline{ span}(x_k)_{k\neq j}, \forall j \in \mathbb N .$$
We say that $(x_k)_{k\geq1}$ has a biorthogonal if there exists a sequence  $(z_k)_{k\geq1}$ in $\mathcal H$ such that  
   $$ (x_m,z_n)= \delta_{mn} \ \ \mbox{(Kronecker's $\delta$ symbol}),$$ and in this case we say that $(x_k)_{k\geq1}$ and $(z_k)_{k\geq1}$ are biorthogonal.

 The next lemma is stated in \cite [Lemma 3.3.1] {Ch}.
\begin{Le}
Let $(x_k)_{k\geq1}$ be a sequence in a Hilbert space $\mathcal H.$ Then 
\begin{enumerate}
\item $(x_k)_{k\geq1}$ has a biorthogonal $(z_k)_{k\geq1}$ if and only if $(x_k)_{k\geq1}$ is minimal.
\item If a biorthogonal sequence for $(x_k)_{k\geq1}$ exists, then it is uniquely determined if and only if $(x_k)_{k\geq1}$ is complete in $\mathcal H.$
\end{enumerate}
\end{Le}

A sequence  $(x_k)_{k\geq1}$ is called a Bessel sequence if there exists a constant $b>0$ such that 
 $$ \sum_{k=1} ^{\infty} |(x,x_k)_{\mathcal H}|^2 \leq  b \ \|x\|^2, \ \forall x\in \mathcal H.$$
 The constant $b$ is called a Bessel bound or an upper bound for $(x_k)_{k\geq1},$  and the smallest upper bound $b$  for $(x_k)_{k\geq1},$  will be denoted $b_X.$
The following  lemma stated in \cite [Theorem 7.4] {H}, characterizes all Bessel sequences for $\mathcal H$ starting with one orthonormal basis.
 \begin{Le}
 Let $(w_k)_{k\geq1}$ be an orthonormal basis for $\mathcal H.$ Then the Bessel sequences for $\mathcal H$ are precisely the sequences $(Uw_k)_{k\geq1},$ where  $U$ is a bounded linear operator on $\mathcal H.$ 

\end{Le}

A  sequence $(x_k)_{k\geq1}$ in a Hilbert space $\mathcal H$  is said to be a Riesz basis for $\mathcal H$
if there exists an orthonormal basis $(w_k)_{k\geq1}$   for $\mathcal H$ and an isomorphism  $U$ on $\mathcal H$ such that 
 $$ \forall k\geq1, \ x_k =Uw_k.$$ 
The next theorem stated in \cite [Theorem 3.6.6] {Ch}, gives equivalent conditions for $(x_k)_{k\geq1}$ being a Riesz basis.
\begin{thm}
For a sequence $(x_k)_{k\geq1}$ in  a Hilbert space $\mathcal H,$ the following statements are equivalent:
\begin{enumerate}
\item $(x_k)_{k\geq1}$ is a Riesz basis for $\mathcal H.$
\item $(x_k)_{k\geq1}$ is complete in $\mathcal H$ and there exist $a,b>0$ such that for all finite scalar sequence $(c_k)_{k\geq1}$
$$ a \sum_{k=1}^{\infty} |c_k| ^2 \leq \sum_{k=1}^{\infty} \| c_k x_k \|^2 \leq b \sum _{k=1}^{\infty} |c_k| ^2 .$$
\item $(x_k)_{k\geq1}$ is a complete Bessel sequence, and has a complete biorthogonal sequence $(y_k)_{k\geq1}$ which is also a Bessel sequence.
\end{enumerate}
\end{thm}

For a given sequence $X=(x_k)_{k\geq1}$ in $\mathcal H,$ let us introduce some related operators. The synthesis operator associated with  $X=(x_k)_{k\geq1}$ is defined as follows:
$$  \mathcal D(S_X)= \{  (c_k)_{k\geq1} \in \ell^2(\mathbb N^*)  \ / \ \displaystyle{\sum_{k}} c_k x_k \ \mbox{converges} \},$$
and for $(c_k)_{k\geq1} \in  \mathcal D(S_X),$ 
$$   S_X(c_k)_{k\geq1} = \sum_{k=1}^\infty c_k x_k.$$
Since the finite sequences  are dense in $\ell ^2(\mathbb N^*)$ and contained in  $\mathcal D(S_X),$ the synthesis operator $S_X$ is densely defined. The analysis operator associated with the sequence $X=(x_k)_{k\geq1}$ is defined by

$$ \mathcal D(A_X)= \{  x\in \mathcal H \ / \ ((x,x_k)_{\mathcal H})_{k\geq 1}\in \ell^2(\mathbb N^*) \},$$
and for $x \in \mathcal D(A_X),$
$$  A_Xx  = ((x,x_k)_{\mathcal H})_{k\geq 1}. $$
The following lemma is stated in \cite [Theorem 7.4] {H}.
\begin{Le}
 Let $(x_k)_{k\geq1}$ be a sequence in $\mathcal H.$ Then, $(x_k)_{k\geq1}$ is a Bessel sequence if and only if the associated synthesis operator  is bounded. 
 \end{Le}
The next lemma is stated in   \cite [Lemma 3.1] {CCLA} and \cite [Lemma 8.4.2] {Ch}.
\begin{Le} 
Let $\mathcal H$ be a Hilbert space and $X=(x_n)_{n\geq 1}$ an arbitrary sequence in $\mathcal H.$ Then, the following hold
\begin{enumerate}
\item The analysis operator $A_X$ is closed.
\item If the analysis operator $A_X$ is densely defined, then the adjoint operator $A_X^*$ is an extension of the synthesis operator $S_X,$ i.e., $ S_X \subset A_X^*.$

\end{enumerate}

\end{Le}

Note that if $(x_k)_{k\geq1}$ is an orthonormal basis  for $\mathcal H,$ the associated analysis operator $A_{X}$ is a unitary isomorphism. 
In the case $(x_k)_{k\geq1}$ is a Bessel sequence,   $A_X$  is bounded and 
$$ S_X= A_X^*,$$ 
where $S_X$ is the synthesis operator associated with $(x_k)_{k\geq1}.$

The rest of this section concerns some basic results for Sobolev spaces in Lipschitz domains. 
\subsection{Sobolev spaces in Lipschitz domains}
    Throughout this section, $\Om$ is an open subset of $\mathbb R^d,$ $d=1,2,3,...$,  $\p \Om$ its boundary and $\overline{\Om}$  its closure. $\mathcal {C}^k(\Om)$ denotes the space of functions mapping $\Om$ into $\mathbb C$ such that all partial derivatives up to order $k$ are continuous, where $k \in \mathbb Z_+ $ and we denote by $\mathcal C^k(Q),$ for $Q$ a closed subset of $\mathbb R^d,$ the space of restrictions to $Q$ of all functions in $\mathcal C^k(\mathbb R^d).$  \\
    Consider the multi-index $ \alpha = (\alpha_1,...\alpha_d) \in \mathbb Z_+^d.$
     We define $\displaystyle |\alpha|= \sum _{k=1}^d \alpha_k.$ For $f \in \mathcal{C}^m(\Om)$ and $|\alpha|\leq m,$ we define
    $$\partial ^{\alpha}f = \frac{\partial^{|\alpha|}f}{\partial x_1^{\alpha_1}...\partial x_d^{\alpha_d} } = \frac{\partial ^{\alpha_1}}{\partial x_1 ^{\alpha_1}}... \frac{\partial ^{\alpha_d}}{\partial x_d ^{\alpha_d}} f.$$
    If $Q \subset \mathbb R^d$ is compact, we may equip $\mathcal C^k(Q)$ with the norm
    $$ \|\varphi\|_{\mathcal C^k(Q)} = \sup _{x\in Q,  |\alpha| \leq k} |(\partial ^{\alpha}\varphi)(x)|.$$
    \\
   We denote by $\mathcal C ^{\infty}(Q)$ for closed $Q\subset \mathbb R ^d$, the intersection of all $\mathcal C^k(Q),$ for $k \in \mathbb Z _+.$ The closure of the set $ \{ x \in \Om \ | \ f(x) \neq 0 \} $  where $f \in \mathcal C(\Om),$ is called the \textbf{support} of $f$ and denoted $supp f.$ 
     A function $f \in \mathcal C^{\infty}(\Om)$ is said to be a test function if $supp f$ is a compact subset of $\Om $ and the set of all test functions on $\Om$ is denoted by $\mathcal{C}_c^{\infty}(\Om).$ For   a sequence  $(\varphi_n)_{n\geq 1}$   in $\mathcal C_c^{\infty}(\Om)$ and $\varphi \in \mathcal C_c^{\infty}(\Om),$ we say that $(\varphi_n)_{n\geq 1} $ converges to $\phi$ in $\mathcal C_c^{\infty}(\Om)$ if there exists a compact $Q\subset \Om$ such that for all $n\geq 1 $ $supp(\varphi_n) \subset Q$ and for all multi-index $\alpha \in \mathbb Z_{+}^d$ the sequence $(\partial ^{\alpha} \varphi_n) _{n\geq 1} $ converges uniformly to $\partial ^{\alpha} \varphi.$ The space $\mathcal C_c^\infty(\Om)$
induced by this convergence  is denoted $\mathscr D(\Om).$
Moreover,
the action of a linear map  $u:\mathscr{D}(\Om) \longrightarrow \mathbb C$ on the test function $\varphi$ is denoted by $<u,\varphi>.$

A distribution on $\Om$ is a linear map $u:\mathscr{D}(\Om) \longrightarrow \mathbb C$ such that for all compact $Q\subset \Om $, there exists $m \in \mathbb Z_+$ and $c \geq 0$ such that
     $$  |<u,\varphi >| \leq c \  \|\varphi\|_{\mathcal C^m(Q)} \ \ \forall \varphi \in \mathscr{ D}(\Om),  $$
     where $m$ and $c$ may depend on $Q.$   We denote by $\mathscr {D'}(\Om)$ the vector space of distributions on $\Om.$  For $u\in \mathscr{ D'}(\Omega)$ a distribution,  one can define its partial derivative with respect to $x_i$ to be the distribution $\frac{\partial u}{\partial x_j},$ specified by $$<\frac{\partial u}{\partial x_j}, \varphi> = -<u,\frac{\partial \varphi}{\partial x_j}> \ \ \ \forall \varphi \in \mathscr{ D}(\Omega).$$
 
 Once the derivative has been defined,  it will be easy to define recursively higher derivatives by induction, i.e $\frac{\partial }{\partial x_i}(\frac{\partial u}{\partial x_j}) .$
 
 We denote by $H^k(\Omega)$  the Sobolev space of all distributions $u$ defined  on $\Omega$ such that all partial derivatives of order at most $k$ lie in $L^2(\Omega), i.e., $ $$ \partial ^{\alpha} u \in L^2(\Omega), \ \ \forall  \ |\alpha | \leq k .$$
$H^k(\Omega)$ equipped with the norm
  \begin{equation*}
 \|u\|_{k,\Omega} = \Big (\sum_{|\alpha| \leq k} \int_{\Omega} |\partial ^{\alpha}u |^2 \ dx \ \Big )^{1/2},
  \end{equation*}
 \\  associated with the inner product
 \begin{equation*}
(u,v)_{k,\Omega} = \sum_{|\alpha | \leq k} \int_{\Omega} \partial ^{\alpha}u \ \overline{\partial^{\alpha}v} \ dx, \ \ \ \forall u,v \in H^k(\Omega)
 \end{equation*}
is a Hilbert space, where $\overline{\partial^{\alpha}v} $ is the conjugate of $\partial^{\alpha}v.$ Sobolev spaces $H^s(\Omega)$ for non-integer $s$ are defined by the real interpolation method (see \cite{Ad}, \cite{Mcl} and \cite{T}).
\begin{Def}
 Let $\Omega$ be an open subset of $\mathbb R^d$ with boundary $\partial \Omega$ and closure $\overline{\Om}. $ We say that $\partial \Omega$ is Lipschitz continuous if for every $x\in \partial \Omega$ there exists a coordinate system $(\widehat{y}, y_d)\in \mathbb R^{d-1}\times\mathbb R$, a  neighborhood $Q_{\delta,\delta'}(x)$  of $x $ and a Lipschitz function $\gamma_x:\widehat{Q}_{\delta} \rightarrow \mathbb R$  with the following properties:
\begin{enumerate}
\item $\Omega \cap Q_{\delta,\delta'}(x) = \{  (\widehat{y},y_d) \in Q_{\delta,\delta'}(x) \ / \ \gamma_x(\widehat{x}) < y_{d} \};$
\item $\partial \Omega \cap  Q_{\delta,\delta'}(x) =\{  (\widehat{y},y_d) \in Q_{\delta,\delta'}(x) \ / \ \gamma_x(\widehat{x}) = y_{d} \};$
\end{enumerate}
where $$ Q_{\delta,\delta'}(x) = \{  (\widehat{y},y_d) \in \mathbb R^d \ / \  \|\widehat{y}-\widehat{x}\| _{\mathbb R^{d-1}} < \delta  \ \ and \ \ |y_d - x_d | < \delta' \ \}$$
and $$ \widehat{Q}_{\delta}(x) = \{ \widehat{y} \in \mathbb R^{d-1} \ / \  \|\widehat{y}-\widehat{x}\| _{\mathbb R^{d-1}} < \delta  \} $$
for  $\delta, \delta' > 0.$ 
An open connected subset $\Omega \subset \mathbb R^d$ whose boundary is Lipschitz continuous is called a Lipschitz domain.
\end{Def}
If $\Om$ is a Lipschitz hypograph, then according to Mclean \cite{Mcl}, we can construct Sobolev spaces on its boundary $\p \Om$ in terms of Sobolev spaces on $\mathbb R^{d-1}
,$ as follows. For $g\in L^2(\p \Om),$ we define
$$ g_{\gamma}(\widehat{x}) = g(\widehat{x}, \gamma(\widehat{x})) \ \mbox{for} \  \widehat{x} \in \mathbb R^{d-1}  ,$$
 put $$H^s(\p \Om) = \{ \ g\in L^2(\p \Om) \ | \ g_{\gamma} \in H^s(\mathbb R^{d-1})  \ \mbox{for} \ 0\leq s\leq 1 \}, $$
and equip this space with the inner product
$$ (g,y)_{s,\p \Om} = (g_{\gamma}, y_{\gamma}) _{s,\mathbb R^{d-1}}.$$
where
$$(u,v)_{s,\mathbb R^{d-1}} = \int_{\mathbb R^{d-1}} (1+|\xi| ^2) ^s \widehat{u}(\xi) \overline{\widehat{v}(\xi)} \ d\xi .$$

Recalling that any Lipschitz function is almost everywhere differentiable so, any Lipschitz hypograph $\Omega$ has a surface measure $\sigma,$ and an outward unit normal $\nu$ that exists $\sigma$-almost everywhere on $\p \Om.$ If $\Om$ is a Lipschitz hypograph then

$$ d\sigma(x) = \sqrt{1+\| \nabla \gamma (\widehat{x}) \| _{\mathbb R^{d-1}}^2  } d \widehat{x} $$
and $$\nu(x) = \frac{  ( - \nabla \gamma (\widehat{x}),1)}  { \sqrt{1+\| \nabla \gamma (\widehat{x}) \| _{\mathbb R^{d-1}}^2  } }$$
for almost every $x\in \p \Om.$ \\
Suppose now that $\Om$ is a Lipschitz domain. Since $\p \Om \subset \bigcup _{x\in \p \Om}  Q_{\delta , \delta'}(x) $ and  that $\p \Om$ is compact, there exist then $x^1, x^2,...,x^n \in \p \Om $ such that $$\p \Om \subset \bigcup _{j=1}^n Q_{\delta , \delta'}(x^j) .$$ It follows that the family $ (W_j) = (Q_{\delta,\delta'}(x^j))$ is a finite open cover of $\partial\Omega,$  i.e., each $W_j$ is an open subset of $\mathbb R^d$, and $\p \Om \subseteq  \bigcup _j W_j .$

 Let $(\varphi_j) $ be a partition of unity subordinate to the open cover $(W_j)$ of $\p \Om,$ i.e.,
$$\varphi_j \in \mathscr D(W_j) \ \ \mbox{and} \ \sum _j  \varphi_j(x) =1  \ \ \mbox{for all} \ x\in \p\Om.$$
The inner product in $H^s(\p \Om)$ is then defined by
$$(u,v)_{H^s(\p \Om)} = \sum_j (\varphi_j u, \varphi_j v) _{H^s(\p \Om_j)},$$
where $\Om_j$ can be transformed to a Lipschitz hypograph by a rigid motion, i.e., by a rotation plus a translation and satisfies
$$W_j \cap \Om = W_j \cap \Om_j \ \mbox{for each}  \ j.$$ 
 It is interesting to mention that a different choice of $(W_j), (\Om_j)$ and $(\varphi_j)$ would yield the same space $H^s(\p \Om)$ with an equivalent norm, for $0\leq s \leq 1.$ For further lectures see (\cite{Ad} and \cite{Mcl}) .

The following lemmas are stated in (\cite{Gr1} and \cite{N}).

\begin{Le}
For a bounded Lipschitz domain $\Omega$  with boundary $\partial \Omega,$ the space $H^{1/2}(\partial \Omega)$ is dense in $L^2(\partial \Omega).$ 
\end{Le}
\begin{Le}  Let $\Omega$ be a bounded Lipschitz domain in $\mathbb R^d.$ Then, the space $H^s(\Omega)$ is compactly imbedded in $H^{s'}(\Omega)$ for all $s' < s$ in $\mathbb R.$
\end{Le} 
\begin{Def}
For a bounded Lipschitz domain with boundary $\partial \Omega,$ the space $H^{-1/2}(\partial \Omega)$ is the dual space of $H^{1/2}(\partial \Omega).$
\end{Def}
Several mathematicians contributed to the study of the trace spaces in Lipschitz domains, most notably Gagliardo  on $W^{1,p}(\Omega)$  for $1\leq p \leq +\infty $ (see \cite{Ga}) and  Costabel on $H^s(\Omega)$ for $\frac{1}{2}< s < \frac{3}{2}$ (see \cite{Co}). \\

Throughout the rest of this paper,  $\Om$  denotes a bounded Lipschitz domain of $\mathbb{R}^{d}$.
\section{The main key ingredients}

Let $\Om$  be a bounded Lipschitz domain of $\mathbb{R}^{d}, d\geq 2$. The trace map takes each  continuous function $u$ on $\overline{\Omega}$ to its restriction on $\p\Om.$ Under the condition $\Om$ is a bounded Lipschitz domain, this trace map may be extended to be a continuous surjective operator denoted $\Gamma_s$  from $H^s(\Om)$ to $H^{s-1/2}(\p\Om),$ for $\frac{1}{2}<s<\frac{3}{2}$  (see \cite{Ad}, \cite{Co}, \cite{Mcl} and \cite{N}).
The range space and the null space of $\Gamma_s$ are respectively given by
	$$\mathcal R(\Gamma_s)= H^{s-1/2}(\partial \Omega) \ \mbox{and} \ \mathcal N(\Gamma_s) = H_0^s(\Omega),$$
	where $H_0^s(\Omega)$ is  the closure in $H^s(\Omega)$ of infinitely differentiable functions compactly supported in $\Omega.$\\

\ \ \ \ \ Let us set  $\Gamma=T_1 \Gamma_1,$  where $\Gamma_1$ is the trace operator from $H^1(\Omega)$ into $H^{1/2}(\partial \Omega)$ and $T_1$ is the embedding operator from  $H^{1/2}(\partial \Omega)$  into $L^2(\partial \Omega).$  According to Gagliardo (see \cite{Ga}), it follows that
$\mathcal R(\Gamma) = H^{1/2}(\p \Om) $ and 
 $\mathcal N(\Gamma)= H_0^1(\Om).$ 
Since $\Gamma_1$ is bounded and $T_1$ is compact (see \cite{N}),  $\Gamma$ is compact. Moreover, since $\mathcal R(\Gamma)$ is dense in $L^2(\p \Om),$ we have the following lemma:
\begin{Le}
Let $\Gamma$ be the trace operator from $H^1(\Om)$ into $L^2(\p \Om).$ Then, the adjoint operator $\Ga^*$ is injective and compact.
\end{Le}

Now, we induce $H^1(\Om)$ by the following inner product

  \begin{equation}
  (u,v)_{\p, \Om}= \int_{\Om} \nabla u \nabla v dx + \int_{\p \Om} \Ga u \Ga v d \sig ~~~~ \forall ~~u,v ~~ \in H^1(\Om).
  \end{equation}

  The associated  norm $\|.\|_{\p, \Om}$ is given by
 \begin{equation}
  \|u\|_{\p, \Om}=\left(\|\nabla u\|^2_{L^2(\Om)}+\|\Ga u\|^2_{L^2(\p\Om)} \right)^{1/2},
   \end{equation}

and  $H^{1}(\Omega)$ induced with the inner product $(.,.)_{\partial,\Om}$ will  be denoted $H_{\partial}^{1}(\Omega).$ 
For $v\in C^1(\overline{\Om})$ the normal derivative map $\partial _{\nu},$ maps each $v$ to $\partial_{\nu} v = \nu. (\nabla v)_{| \p\Om}$ onto $L^2(\p\Om).$  Moreover, under the condition $\Omega$ is a bounded Lipschitz domain, $\partial_{\nu}$  may be extended to be a bounded linear operator denoted $\widehat \partial _{\nu} $ from $H_{\Delta}^1(\Om)$ to $H^{-1/2} (\p\Om)$ (see \cite{Gr1}).
In the following, we recall Green's formula (see \cite{Gr1} and \cite{N}).
\begin{pr} \textbf{"Green's formula"}  
Let $\Om$ be a bounded Lipschitz domain, then for all $u \in H_{\Delta}^1(\Om)$ and $ v \in H^1(\Om)$ one has
$$ \int_{\Om} \nabla u \nabla v dx = - \int_{\Om} \Delta u \ Ev  \ dx + <\widehat{\partial _{\nu}} u, \Gamma_1 v> ,$$
where $E$ is the embedding operator from $H^1(\Om)$ into $L^2(\Om)$ and $<.,.>$ is the duality pairing between $H^{-1/2}(\p \Omega)$ and $H^{1/2}(\p \Omega).$
\end{pr}

The following proposition characterizes $\Gamma^*.$
\begin{pr}

 For all $g \in L^2(\p \Om)$,  $\Gamma^*$ is the  solution operator  of Robin problem for the following Laplace equation
 \begin{equation*}
 \begin{cases}
   \Delta z =0 & \text{  }  (\Om) \\
   \p_{\nu}z+ \Ga z= g & \text{}  (\p \Om),
 \end{cases}
 \end{equation*}
 where $\partial_{\nu}$ is the normal derivative operator considered as non-bounded from $H_{\Delta}^1(\Omega)$ to $L^2(\partial \Omega).$
  \end{pr}
\begin{pre}
 Let  $g \in L^2(\p \Om)$ and  $z=\Ga^*g.$ We have:
    \begin{eqnarray*}
 (\Ga^*g, v)_{ \p ,\Om} & = & \int_{\Om} \nabla z \nabla v dx + \int_{\p \Om} \Ga z \Ga v d \sig  \ \ \ \ (*) \\
&=&  \int_{\p \Om} g \ \Ga v \  d \sig
  \end{eqnarray*}  \\ so that if $ v \in  H^1_0(\Om) ={\cal N}(\Ga)$, then we obtain
  $$  \int_{ \Om} \nabla z \nabla v \ dx = 0  .$$
Since the previous equality characterizes the $H^1-$harmonic functions, then we may write:
 $$ \Delta z = 0 \mbox{   in } \mathscr{D'}(\Omega).$$

Applying Green's formula to $(*),$  we obtain that

\begin{eqnarray*}
\int_{\Om} \nabla z \nabla v dx + \int_{\p \Om} \Ga z \Ga v d \sig
 &=& <\widehat{ \p }_{\nu}z, \Gamma_1 v>  + \int_{\p \Om} \Gamma z \Gamma v d \sig \\
&=&  \int_{\p \Om} g  \Ga v d \sig ,\\
 \end{eqnarray*}
which leads to the following duality pairing on  $H^{1/2}(\p \Omega) \times H^{-1/2}(\p \Omega) $ $$ <\widehat{\partial}_{\nu}z+  \widehat{\Gamma z}, \Gamma_1 v> =<\widehat{g},\Gamma _1 v>  ,$$
where $\widehat{y}$ denotes the embedding of an element $y\in L^2(\p \Om)$ in $H^{-1/2}(\p \Omega).$

 Viewing
 $\mathcal R (\Gamma_1) = H^{1/2}(\partial \Omega)$, it follows that
 $$\widehat{ \p} _{\nu}z + \widehat{\Gamma z} =\widehat{g}, $$ so
 $$\widehat{ \p} _{\nu}z = \widehat{g-\Gamma z}.$$
 Consequently,
$ \widehat{\p}_{\nu} z $ belongs to the range of the embedding operator from $L^2(\p \Omega)$ into $H^{-1/2}(\p \Omega),$
 which means that $\p _{\nu}z \in L^2(\p \Om)$ and that
$$ \partial_{\nu}z+ \Gamma z =g  .\ \  \square$$ 
\end{pre}\\

The trace operator $\Gamma$ being bounded, one considers its Moore-Penrose inverse which we denote by $\Lambda = \Ga^ \dagger \in {\cal C}( L^2(\p\Om), H_{\p}^1(\Om)),$  such that $$ {\cal D}(\Lambda)= {\cal R}(\Ga)=H^{1/2}(\p \Om)  \ \mbox{and} \       {\cal N}(\Lambda^*)= \mathcal N(\Ga)= H_0^1(\Om).$$
Moreover, $\La$ is characterized by the following.
 \begin{pr}
Let $\Ga$ be the trace operator from $H_{\p} ^1(\Om)$ into $L^2(\p \Om)$ and $\La$ its Moore-Penrose inverse. Then, $\Lambda$ is the solution operator of the Dirichlet problem for the Laplace equation with data in $ H^{1/2}(\p \Om).$ Moreover, we have:
$$
   \mathcal R(\Lambda)= \mathcal H^1(\Omega),
   $$
   where
  $\mathcal H^1(\Omega)= \{   v\in H^{1}(\Om) ~~~/~~ \Delta v =0 \ \ \mbox{in} \ \ \mathscr D '(\Om) \}.
$
\end{pr}
 \begin{pre}  Since $\Gamma$ is bounded, it follows that its Moore-Penrose inverse $ \Lambda $ is closed and densely defined with closed range. Moreover, from Lemma 3.1, $\Gamma^*$ is injective, which implies that $\mathcal D(\Lambda)= \mathcal R(\Ga).$
 Also, for $g \in \mathcal D(\Lambda)$ let $v=\Lambda g.$ For $w \in {\cal D}(\Lambda^*),$ we have
 \begin{equation*}
 (v,w)_{\partial,\Omega}= \int_{\Om} \nabla v \nabla w dx + \int_{\p \Om} \Ga v \Ga w d\sigma = (\Lambda g,w)_{\partial,\Omega}=\int_{\partial\Omega}g\Lambda^* w d\sigma,
 \end{equation*}
 so that if $w\in {\cal N}(\Lambda^*)= \mathcal N(\Ga) =H^1_0(\Omega),$ the following holds
 \begin{equation*}
 \int_\Omega \nabla v\nabla w dx =0.
 \end{equation*}
Since the previous equality holds for all $w \in H^1_0(\Omega)$ and  characterizes the $H^1-\mbox{harmonic functions}$, it follows that
  \begin{equation*}
        \begin{cases}
          \Delta v=0 & \text{  }  (\Om) \\
          \Gamma v =g & \text{}  (\p \Om). \ \ \ \ \ \square
        \end{cases}
        \end{equation*}
        \end{pre} 

We now consider the embedding operator:
 $$
  \begin{array}{ccccc}
  E &: & H_{\partial}^1(\Om) & \to & L^2(\Om) \\
   & & v & \mapsto & Ev \\
  \end{array}$$
  which maps each $v\in H^1(\Omega)$ to itself into $L^2(\Omega)$ but obviously with different topologies.
 $H^1(\Om)$ is induced with the inner product $(.,.)_{\p,\Om}$ and $L^2(\Omega)$ with its usual inner product $(.,.)_{0,\Om}$.
The following theorem characterizes $E^*.$
\begin{thm} Let $\Om$ be a bounded Lipschitz domain of $\mathbb R^d$ and $E$  the embedding operator  from $H_{\partial}^1 (\Omega)$ into $L^2(\Omega).$ Then for  $ f\in L^2(\Om), $ the adjoint operator $E^*$ is the solution operator of Robin problem for the following Poisson equation
   \begin{equation}
   \begin{cases}
     -\Delta u =f & \text{  }  (\Om) \\
     \p_{\nu}u+ \Ga u=0& \text{}  (\p \Om).
   \end{cases}
   \end{equation}

\end{thm}
\begin{pre}
Let $ f\in L^2(\Om)$ and $v\in H^1(\Om).$  Putting  $u=E^*f,$  one has
\begin{equation}
  \int_{\Omega} f Ev dx = (E^*f,v)_{\p,\Om}
   = \int_{\Om} \nabla u \nabla v \  dx + \int_{\p \Om} \Ga u \ \Ga v \ d \sigma.
\end{equation}
  Now, if $v \in \mathcal{ C}_ c^{\infty}(\Omega)$ then,
   \begin{eqnarray*}
   (E^*f,v)_{\p,\Om} &= & \int_{\Om}f Ev dx\\
&=&  \int_{\Om} \nabla v \nabla u dx\\
   &=& < -\Delta u, v >_{\mathscr{ D'}(\Omega),\mathscr{ D}(\Omega)}.
    \end{eqnarray*}
  Therefore, $$ f = -\Delta u  \ \ \ \ \mbox{in} \ \ \ \mathscr D'(\Omega).$$

  Applaying Green's formula to (9), one has

   \begin{eqnarray*}
 \int_{\Omega} f \ Ev  dx &=& -\int_{\Om}  Ev \ \Delta u  \ dx + < \widehat{\partial_{\nu}} u , \Gamma_1  v > + \int_{\p \Om} \Ga v \Ga u \ d \sig \\
 &=&   \int_{\Omega} f \ Ev  dx +  < \widehat{\partial_{\nu}} u + \widehat{\Gamma u} , \Gamma_1 v >,
 \end{eqnarray*}
 so far,

$$ < \widehat{\partial_{\nu}} u + \widehat{\Gamma u} , \Gamma_1 v > =0  \ \ \forall v \in H^1(\Omega).$$

Moreover, since $\mathcal R(\Gamma_1)= H^{1/2}(\p \Om),$ it follows that
  \begin{eqnarray*}
\widehat{\partial_{\nu}} u + \widehat{\Gamma u}   =0 \ \ \  \mbox{in} \ \  H^{-1/2} (\partial \Omega).
    \end{eqnarray*}
Consequently, $\widehat{\p _{\nu}} u $ belongs to the range of the embedding operator acting from $L^2(\p \Omega)$ to $H^{-1/2}(\p \Omega)$ and $\partial _{\nu} u \in L^2(\p \Omega),$ which implies that $$ \partial_{\nu} u + \Gamma u  =0 \ \mbox{ in} \   L^2(\partial \Omega) .  \ \ \ \square $$
\end{pre}

Now, for $f\in L^2(\Omega),$ let us consider $E_0^*$ the solution operator of the Dirichlet problem for the following Poisson equation
 \begin{equation*}
              \begin{cases}
                -\Delta u^0=f & \text{  }  (\Omega) \\
                \Gamma u^0=0 & \text{}  (\partial \Omega)
              \end{cases}.
              \end{equation*}

 By setting $E_1^*= E^*-E_0^*$ and $u^1 =E_1^*f,$ it follows that $u^1$ is a solution of  the Dirichlet problem for the  following Laplace equation
 \begin{equation*}
              \begin{cases}
                -\Delta u^1=0 & \text{  }  (\Omega) \\
               \Gamma u^1=\Gamma u & \text{}  (\partial \Omega),
              \end{cases}
              \end{equation*}
   where  $u$ is the solution of (8). Furthermore,
using  Rellich-Ne\v{c}as Lemma, we can prove the following crucial theorem:

 \begin{thm} Let $\Gamma$ be the trace operator from $H_{\p}^1(\Om)$ into $L^2(\p \Om)$ and $E_1^*$ defined as above. Then,  we have
 $  E_1^* = \Gamma ^* K^*,$
 where $K$ is the  solution operator of the Dirichlet problem for the Laplace equation (1).
 \end{thm}
 
 \begin{pre} Putting  $u^1 = E_1^*f,$ we have $u^1= u- u^0,$ where $u^0$ and $u$ are solutions of the followings problems
 \begin{equation*}
               \begin{cases}
                 -\Delta u^0=f & \text{  }  (\Omega) \\
                 \Gamma u^0=0 & \text{}  (\partial \Omega)
               \end{cases}
               \end{equation*}
 and   \begin{equation*}
    \begin{cases}
      -\Delta u =f & \text{  }  (\Om) \\
      \p_{\nu}u+ \Ga u=0& \text{}  (\p \Om)
    \end{cases}
    \end{equation*}
respectively.
 Since $u^1,u^0 \in H_{\Delta}^1(\Omega),$  it follows that
 $$\widehat{\partial}_{\nu} u^1, \widehat{\partial}_{\nu} u^0 \in H^{-1/2}(\partial \Omega),$$
 therefore, $$\widehat{\partial}_{\nu}u = \widehat{\partial}_{\nu} u^0 +\widehat{\partial}_{\nu} u^1 \in  H^{-1/2}(\partial \Omega).$$ On the one hand,  $\partial_{\nu} u = -\Ga u$ implies that $\partial_{\nu} u \in L^2(\p \Om) ,$  and by Rellich-Ne\v{c}as Lemma, we have  $\partial_{\nu} u^0 \in L^2(\p \Om) .$ Which implies  that $\partial_{\nu} u^1 \in L^2(\p \Om).$
On the other hand, since $$  \partial_{\nu} u^1  +\Gamma u^1 =  \partial_{\nu} u^1   -  \partial_{\nu} u=- \partial_{\nu} u^0 ,$$ 
and that the adjoint operator $K^*$   takes each $f \in L^2(\p \Omega)$ to $-\partial_{\nu} u^0$ onto $L^2(\p \Omega),$ where $u^0$ is the solution of the Dirichlet problem for the Poisson equation $(3),$ it follows that 
$$ \partial_{\nu} u^1  +\Gamma u^1 = K^* f ,$$ and that $u^1$ is  the unique solution of 
\begin{equation*}
   \begin{cases}
     -\Delta u^1 =f & \text{  }  (\Om) \\
     \p_{\nu}u^1+ \Ga u^1= - \partial_{\nu} u^0 & \text{}  (\p \Om).
   \end{cases}
   \end{equation*}

Therefore,  $$E_1^*= \Ga^* K^*. \ \ \ \ \ \ \ \ \ \ \ \ \ \ \ \ \ \  \square$$ 
 \end{pre} \\

Now, the operator $E_1$ being bounded, one considers its Moore-Penrose inverse which we denote by $F_1$ such that $$\mathcal D(F_1) = \mathcal R(E_1) \oplus \mathcal N(E_1^*),$$

and 
  $$\mathcal R(F_1)=\mathcal H^1(\Omega) \ \ \mbox{and} \ \ \mathcal N(F_1)= \mathcal N(E_1^*)=H_0^1(\Omega).$$

 According to Proposition 2.3, 
 the operator $ F_1^* (I+F_1 F_1^*)^{-1/2}$ acting from $H^1(\Om)$ into $L^2(\Om)$  is bounded with closed range, i.e., 
$$\mathcal R  (F_1^* (I+F_1 F_1^*)^{-1/2}) = \mathcal R (F_1^*)= \mathcal H(\Om)$$
and
$$\mathcal N( F_1^* (I+F_1 F_1^*)^{-1/2}) = \mathcal N (F_1^*) = H_0^1(\Om) ,$$
and from Corollary 2.1, we have the following lemma.
 \begin{Le} The operator
    $ F_1^* (I+F_1 F_1^*)^{-1/2}$ is an isomorphism from  $\mathcal H^1(\Om)$ into $\mathcal H(\Om).$
 \end{Le} 
Let us now set
$$\Gamma_0^* = F_1^* (I+F_1 F_1^*)^{-1/2} \Gamma ^* .$$ The following lemma characterizes $\Ga_0^*.$

\begin{Le}
The operator $\Gamma_0^*$ defined above  is compact and injective.
\end{Le}
\begin{pre}
Knowing that $\mathcal R  (F_1^* (I+F_1 F_1^*)^{-1/2}) = \mathcal H(\Om),$
we have $ \mathcal R(\Gamma_0^*) \subset \mathcal H(\Omega). $ Moreover, $\Ga$ being compact (see $\S 3$), it follows  by Schauder's theorem (Theorem 2.1) that $\Ga^*$ is compact as well. Therefore, the boundedness of $F_1^* (I+F_1 F_1^*)^{-1/2} $ implies  that  $\Ga_0^*$ is compact. The injectivity of $\Gamma_0^*$ holds for the reason that $\Gamma^*$ is injective and that $\mathcal R(\Gamma^*)\subset \mathcal N(F_1^* (I+F_1 F_1^*)^{-1/2})^{\perp}.$ $\square$
\end{pre}
\\ \\
 \ \ \ \ \ \ \ \  Let us now return to the Dirichlet problem for the Laplace equation (1), where we have considered its solution operator $K$ and its adjoint $K^*.$ 
Composing  $\Gamma_0^*$ by $K^*,$ we obtain
  \begin{eqnarray*}
  \Gamma_0^*  K^* &= & F_1^* (I+F_1 F_1^*)^{-1/2} \Gamma ^* K^* , 
  \end{eqnarray*}
and in view of  Theorem 3.2, we have $E_1^* = \Gamma^* K^*,$ which leads to
 \begin{equation*}
   \Gamma_0^*  K^* = F_1^* (I+F_1 F_1^*)^{-1/2} E_1^*  .
   \end{equation*}
   On the other hand, since $(I+F_1 ^* F_1)^{-1/2} F_1^* \subset  F_1^* (I+F_1 F_1^*)^{-1/2}$ and that  $\mathcal R(E_1^*)\subset \mathcal D(F_1^*),$ one has
  $$
   \Gamma_0^*  K^*= (I+F_1 ^* F_1)^{-1/2} F_1^* E_1^* , $$
   and viewing $F_1^* E_1^* =P_{\mathcal H(\Omega)},$ it follows that
   $$
 \Gamma_0^*  K^*=(I+F_1 F_1^*)^{-1/2} P_{\mathcal H(\Omega)}.$$

  Therefore, $$ \mathcal R(\Gamma _0^* K^*) = \mathcal R((I+F_1^* F_1) ^{-1/2}  P_{\mathcal H(\Omega)}) .$$

  Viewing  $\Gamma_0^*$ is compact,  $K\in \mathcal B(L^2(\partial \Omega),L^2(\Omega))$   and that 
 $ (I+F_1^* F_1) ^{-1/2}  P_{\mathcal H(\Omega)} $ is  self-adjoint, it follows that
 $\Gamma_0 ^* K^*$  is compact  and self-adjoint. Therefore, there exists a sequence $ ((\kappa_n, \phi_n)) _{n\geq 1}$ in  $\mathbb R _+^* \times \mathcal H(\Omega) $ such that for all $n \geq 1,$
  $$ \Gamma _0^* K^* \phi _n =  \kappa _n^2 \phi_n .$$
  Moreover, the sequence $(\phi_n)_{n\geq 1} $ is an orthonormal basis for the Bergman space $\mathcal H(\Omega).$ \\ \\
  The aim of the next section is to prove the main result of this paper (Theorem 1.1).

\section{The sequences $(g_n)_{n\geq 1}$ and $(y_n)_{n\geq 1}$}

As we stated in the first section of this paper, the sequences $(g_n)_{n\geq 1}$ and $(y_n)_{n\geq 1}$ are defined  for all $n\geq 1,$ by
   $$ \Gamma_0^* \phi_n = \kappa_n y_n \ \mbox{and} \ K^* \phi_n = \kappa_n g_n,$$ 
    where 
 $ ((\kappa_n, \phi_n)) _{n\geq 1}$ in  $\mathbb R _+^* \times \mathcal H(\Omega) $ is a sequence of couples associated to $\Gamma_0 ^* K^*.$  A first  remark is that $\Gamma_0^*$ and $K$ satisfie the following
  $$ \Gamma_0^* g_n = \kappa_n \phi_n = K y_n .$$ 
 
Denote by $\mathcal G(\p \Om)$ and $\mathcal Y(\p \Om)$ the closures of $ span (g_n)_{n\geq 1}$ and $span(y_n)_{n\geq 1}$  respectively. The principal objective of this section will be to prove that the sequences $(g_n)_n$ and $(y_n)_n$ are Riesz bases. This will be the key ingredient to prove the $H^{1/2}$ regularity result for the problem(1), stated in Theorem 1.2.
\begin{Le}
The sequences $(y_n)_{n\geq 1}$ and $(g_n)_{n\geq 1}$ are biorthogonal and $(y_n)_{n\geq 1}$ is complete.
\end{Le}
\begin{pre}
Let $m,n \geq 1.$ We have 
\begin{eqnarray*}
 (g_n, y_m)_{0,\p \Omega} &=& \frac{1}{\kappa_n \kappa_m}(\kappa_n g_n, \kappa_m y_m)_{0,\p \Omega}  \\
 &=& \frac{1}{\kappa_n \kappa_m}(K^* \phi_n, \Gamma_0 \phi_m )_{0,\p \Omega} \\
  &=& \frac{1}{\kappa_n \kappa_m}(\Gamma_0^* K^* \phi_n,  \phi_m )_{0, \Omega} \\
 &=& \frac{\kappa_n^2}{ \kappa_n \kappa_m}( \phi_n, \phi_m )_{0,\Omega}.
 \end{eqnarray*}
Since $(\phi_n)_{n\geq 1}$ is an orthonormal basis for $\mathcal H(\Om),$ it follows that $(g_n, y_m)_{0,\p \Omega}= \delta_{nm},$  therefore $(g_n)_{n\geq 1}$ and $(y_n)_{n\geq 1}$ are biorthogonal. To prove that $(y_n)_{n\geq 1}$ is complete, one standard way is to consider an element $g\in L^2(\p \Om)$ such that for all $n\geq 1,$ 
\begin{equation}
(g, y_n)_{0,\p \Om} =0
\end{equation}
and prove that $g=0.$  Multiplying $(10)$ by $\kappa_n,$ it follows that
$$ 0= \kappa_n (g,y_n)_{0,\p \Omega}= (g,\kappa_n y_n)_{0,\p \Omega}= (g,\Gamma_0 \phi_n)_{0,\p \Omega}=(\Gamma_0 ^*g, \phi_n)_{0,\Omega}$$
therefore, $(\Gamma_0 ^*g, \phi_n)_{0,\Omega}=0$ for all $n\geq 1,$ and since $(\phi_n)_{n\geq 1}$ is an orthonormal basis for $\mathcal H(\Om)$ and that $\Gamma_0 ^*g \in \mathcal H(\Om),$ we obtain that $\Gamma_0 ^*g=0$ which implies that $g=0$ according to the injectivity of $\Gamma_0^*$ from Lemma 3.3. \ $\square$
\end{pre} \\
 
The following corollary is a consequence of Lemma 4.1 and Lemma 2.2.
\begin{co}
The sequences  $(g_n)_{n\geq 1}$ and $(y_n)_{n\geq 1}$ are minimal. 
\end{co}

\begin{Le}
$\mathcal R(K^*)$ and $\mathcal R(\Gamma_0)$ are densely imbedded in $\mathcal G(\p \Omega)$ and $ \mathcal Y(\p \Om)$ respectively. Moreover, $$\mathcal Y(\p \Om)= L^2(\p \Om). $$ 
\end{Le}
\begin{pre}
First, let us prove that $\mathcal R(K^*) \subset \mathcal G(\p \Omega).$ Since $(\phi_n)_{n\geq 1}$ is an orthonormal basis for $\mathcal H(\Omega),$ we have for all $v\in \mathcal H(\Omega)$
$$v=\sum_{n=1}^{\infty} (v,\phi_n)_{0,\Omega} \ \phi_n = \lim_{n \rightarrow {+\infty}} \sum_{j=1}^{n} (v,\phi_j)_{0,\Omega} \ \phi_j,$$
and since $K^*$ is bounded, it follows that
\begin{eqnarray*}
K^*v &=& \lim_{n \rightarrow {+\infty}} \sum_{j=1}^{n} (v,\phi_j)_{0,\Omega} \ K^* \phi_j \\
&=& \lim_{n \rightarrow {+\infty}} \sum_{j=1}^{n} (v,\phi_j)_{0,\Omega} \ \kappa_j g_j.
\end{eqnarray*}
Therefore, $K^*v$ is the limit of a linear combination sequence of elements of $(g_n)_{n\geq 1}.$ Hence $K^*v \in \mathcal G(\p \Omega),$ then  $\mathcal R(K^*) \subset \mathcal G(\p \Omega)$  and since $span (g_n)_{n\geq 1} \subset \mathcal R(K^*),$  it follows that $\mathcal R(K^*)$ is densely imbedded in $\mathcal G(\p \Om).$ In a similar way, we obtain that $\mathcal R(\Gamma_0)\subset \mathcal Y(\p \Omega).$ Moreover, here is another way to obtain the completeness of the sequence $(g_n)_{n\geq 1}:$ since the operator $\Gamma_0^*$ is injective according to Lemma 3.3, $\mathcal R(\Gamma_0)$ is then  dense in $L^2(\p \Om).$ Therefore, we obtain that $\mathcal Y(\p \Om) = L^2(\p \Om).\square$ 
\end{pre}

\begin{Le}
The sequences $(\kappa_n g_n)_{n\geq 1}$ and $(\kappa_n y_n)_{n\geq 1}$ are Bessel sequences.
\end{Le}
\begin{pre} Viewing $(\phi_n)_n$ is an orthonormal basis for $\mathcal H(\Omega)$ and that the operators $K^*$ and $\Gamma_0^*$ are bounded, it follows by Lemma 2.3 that  $(\kappa_n g_n)_{n\geq 1}$ and $(\kappa_n y_n)_{n\geq 1}$ are Bessel sequences. $\square$
\end{pre}

\begin{co}
The synthesis  operators associated with the sequences $(\kappa_n g_n)_{n\geq 1}$ and  $(\kappa_n y_n)_{n\geq 1}$ are  bounded.
\end{co}

\begin{co}
The analysis operators associated with the sequences $(\kappa_n g_n)_{n\geq 1}$ and  $(\kappa_n y_n)_{n\geq 1}$ are  bounded.
\end{co}

In the rest of this paper, we denote by
$(A_G,\mathcal D(A_G))$ and $ (A_Y,\mathcal D(A_Y))$   the analysis operators associated with the sequences $(g_n)_{n\geq 1} , (y_n)_{n\geq 1}, $ and by $(S_G,\mathcal D(S_G))$ and $(S_Y,\mathcal D(S_Y))$ their associated synthesis operators respectively. Denote also by  $A_{\Phi}$ the analysis operator associated with the orthonormal basis $(\phi_n)_{n\geq 1}$ and by $M_{\kappa}$ the multiplication operator on $\ell^2(\mathbb N^*)$ 
by the sequence $\kappa=(\kappa_n)_{n\geq 1}$ such that

$$M_{\kappa} : \ell^2(\mathbb N^*) \longrightarrow  \ell^2(\mathbb N^*) ,$$ such that for all scalar sequence $(x_n)_{n\geq 1}\in \ell^2(\mathbb N^*),$ $$
 M_{\kappa} (x_n)_{n\geq 1} = (\kappa_n x_n)_{n\geq1}.$$
 
In particular, for the sequences $G= (g_n)_{n\geq 1} $ and $Y= (y_n)_{n\geq 1} ,$ we adopt the following notation:
 $$M_{\kappa} G =(\kappa_n g_n)_{n\geq1} = \kappa G\ \mbox{and} \  M_{\kappa} Y =(\kappa_n y_n)_{n\geq1} = \kappa Y.$$
 The next lemma will prove to be crucial.
\begin{Le} We have
$$A_Y K^* =  M_{\kappa} A_{\Phi}=A_G \Gamma_0.$$
\end{Le}
\begin{pre}
For $v\in \mathcal H(\Omega),$  we have
\begin{eqnarray*}
M_\kappa A_{\Phi} v &=&  (\kappa_n(v,\phi_n)_{0,\Omega})_{n\geq1}\\
&=& ((v,\kappa_n\phi_n)_{0,\Omega})_{n\geq1}\\
&=& ((v, K y_n)_{0,\Omega})_{n\geq1}\\
&=& ((K^*v,  y_n)_{0,\p \Omega})_{n\geq1}\\
&=& A_Y K^* v.
\end{eqnarray*}
In a similar way, we obtain
\begin{eqnarray*}
M_\kappa A_{\Phi} v &=&  (\kappa_n(v,\phi_n)_{0,\Omega})_{n\geq1}\\
&=& ((v,\kappa_n\phi_n)_{0,\Omega})_{n\geq1}\\
&=& ((v, \Gamma_0^* g_n)_{0,\Omega})_{n\geq1}\\
&=& ((\Ga_0v,  g_n)_{0,\p \Omega})_{n\geq1} \\
&=& A_G\Gamma_0 v. \ \square 
\end{eqnarray*}

\end{pre}
\begin{co}
The analysis operators $(A_G,\mathcal D(A_G)) $ associated  with the sequence $(g_n)_{n\geq1}$ is closed and densely defined.
\end{co}
\begin{pre}
According to Lemma 2.5, $(A_G,\mathcal D(A_G))$ is closed. On the other hand, we have shown in Lemma 4.4, that $\mathcal R(\Ga_0) \subset \mathcal D(A_G),$  and since $\mathcal R(\Ga_0)$ is dense in  $L^2(\p \Om),$  it follows then that the analysis operators $(A_G,\mathcal D(A_G)) $  is  densely defined. $\square$
\end{pre}\\

 We will extensively make use of the following lemma:
\begin{Le}
The following hold
\begin{enumerate}
\item $M_\kappa A_G\subset A_\Phi \K= A_{\kappa G}$
\item $M_\kappa A_Y\subset A_\Phi \Gamma_0^*=A_{\kappa Y}$
\item $K^*\A^*_\Phi=\A_G^*M_\kappa$
\item $\Gamma_0 A^*_\Phi=A_{\kappa Y}^*$

\end{enumerate}
\end{Le}
\begin{pre} Let $g\in {\cal D}(A_G).$  We have
$$
M_\kappa A_G g=(\kappa_n(g,g_n)_{0,\partial\Omega})_{n\geq1}.
$$
Since for all $n\geq1$, 
\begin{eqnarray*}
\kappa_n(g,g_n)_{0,\partial\Omega}&=&(g,\kappa_ng_n)_{0,\partial\Omega}\\
&=&(g,K^*\phi_n)_{0,\partial\Omega}\\
&=&(K g,\phi_n)_{0,\partial \Omega},
\end{eqnarray*}
it follows that for all $g\in {\mathcal D}( A_G),$ $$M_\kappa A_Gg= A_\Phi K g= A_{\kappa G}g$$ and the first inclusion holds. In a similar way, one can prove (2). Having established in Corollary 4.4 that  the analysis operator $(A_G,\mathcal D(A_G)) $ is densely defined and that $A_{\kappa Y}$ is bounded  in Corollary 4.3, one  considers their adjoint operators  $A_G^*, A_{\kappa Y}^*$ respectively. Therefore, the items (3) and (4) hold by considering the adjoints in $(1)$ and $(2).\ \square$ 
\end{pre} \\

The following result will prove useful in the rest of this text.
\begin{Le}
 The operator $A_{\kappa Y} $ is injective and has a dense range.
\end{Le}
\begin{pre}
For $g\in L^2(\p \Omega),$ the equality 
$A_{\kappa Y} g=0 $ implies that for all $n\geq 1, (g, \kappa_n y_n)_{\p \Omega} =0.$ Moreover, viewing $\kappa_n \neq 0,$ it follows that  $ (g,  y_n) _{\p \Omega} =0.$ Since the sequence $Y= (y_n)_{n\geq 1}$ is complete, we get $g=0.$ Therefore, $A_{\kappa Y} $ is injective. On the other hand,  from the second item of  Lemma 4.5, we have  $A_{\kappa Y} = A_{\Phi} \Ga_0^*,$ and since  $ A_{\Phi}$ is a unitary  isomorphism  and that  $\mathcal R(\Ga_0^*)$ is dense in $ \mathcal H(\Om)$, we deduce that $\mathcal R(A_{\kappa Y})$ is dense in $\ell^2(\mathbb N^*). \ \square$
\end{pre}
\begin{co}  The operator $A_{\kappa Y}^* $ is  bounded, injective and has a dense range.
\end{co}
\begin{Le}
 The operator  $(M_\frac{1}{\kappa} A_G, \mathcal D(M_\frac{1}{\kappa} A_G)) $ is the inverse of $A_{\kappa Y}^*$. Moreover, it is surjective and has a surjective adjoint.
\end{Le}
\begin{pre} From Lemma 4.5, we have  $\Gamma_0 A^*_\Phi=A_{\kappa Y}^*,$ which implies that for all $(c_n)_{n\geq1} \in \ell^2(\mathbb N^*),$  
$$ M_{\frac{1}{\kappa}} A_G A_{\kappa Y}^* (c_n)_n = M_{\frac{1}{\kappa}} A_G \Gamma_0 A_{\Phi} ^* (c_n)_n .$$

In view of Lemma 4.4, we have  $A_G \Ga_0 = M_{\kappa} A_{\Phi},$  which leads to
$$  M_{\frac{1}{\kappa}} A_G A_{\kappa Y}^* (c_n)_n = M_{\frac{1}{\kappa}} M_{\kappa} A_{\Phi} A_{\Phi} ^* (c_n)_n .$$
Moreover, since $A_{\Phi}$ is a unitary isomorphism, we have   $$A_{\Phi} A_{\Phi}^* = I_{\ell^2(\mathbb N^*)}, $$ therefore,
$$  M_{\frac{1}{\kappa}} A_G A_{\kappa Y}^* (c_n)_n =(c_n)_n. $$ Hence, we obtain that
 $$\mathcal R(A_{\kappa Y}^*) \subset \mathcal D(M_{\frac{1}{\kappa}}  A_G)$$
 and that $$ M_{\frac{1}{\kappa}} A_G A_{\kappa Y}^*= I_{\ell^2(\mathbb N^*)}.$$
Having previously established in Corollary 4.5 that $A_{\kappa Y} ^*$ is bounded injective with dense range, it follows that $(M_\frac{1}{\kappa}A_G,\mathcal D(M_\frac{1}{\kappa}A_G))$ is its unique inverse. Moreover,  it  is surjective and has a surjective adjoint. \ $\square$

\end{pre}
\begin{co}
The operator $(A_G^* M_\frac{1}{\kappa},\mathcal D(A_G^* M_\frac{1}{\kappa})) $ is the inverse of $A_{\kappa Y}.$
\end{co}
\begin{co} The operators $ (A_G^*,\mathcal D(A_G^*))$ and $ (A_Y^*,\mathcal D(A_Y^*))$ are surjective. Moreover, the sequence $(g_n)_{n\geq 1}$ is complete.
\end{co}
\begin{pre}
Since the operator $(A_G^* M_\frac{1}{\kappa},\mathcal D(A_G^* M_\frac{1}{\kappa}))$ is  surjective, clearly $(A_G^*,\mathcal D(A_G^*))$ is also surjective as well, then $\mathcal R(A_G^*) = L^2(\p \Omega),$ and this implies that the sequence $(g_n)_{n\geq 1}$ is complete, i.e, $$ \mathcal G(\p \Omega) = L^2(\p \Omega).$$ On the other hand, we have $$ span(g_n)_{n\geq 1} \subset \mathcal R(K^*) ,$$ 
and according to Lemma 4.2, $\mathcal R(K^*)$ is densely imbedded in $ \mathcal G(\p \Omega).$
 Therefore, we deduce that 
$ \mathcal R(K^*)$ is dense in $L^2(\p \Omega).$  Moreover, by Lemma 4.4, we have $$ \mathcal R(K^*) \subset \mathcal D(A_Y).$$ Therefore, the operator $(A_Y,\mathcal D(A_Y))$ is densely defined. In analogue with the proof of Lemma 4.7,  instead of  $A_{\kappa Y}^* =A_Y^* M_{\kappa},$  we consider the operator $A_G^* M_{\kappa}=A_{\kappa G}^* ,$   then one can prove that $M_{\frac{1}{\kappa}} A_Y$ is the inverse of  $A_G^* M_{\kappa} $ and this establishes that $(A_Y^*,\mathcal D(A_Y^*))$ is surjective. \ $\square$
\end{pre}

\begin{co}
The solution operator of the Dirichlet problem for the Laplace equation (1) is injective.
\end{co}
The next proposition is an essential step towards the main result of this paper.
\begin{pr} The following hold
$$\Ga_0= S_Y A_YK^*= S_Y A_G \Gamma_0,$$
and $$K^* = S_G A_G \Gamma_0 = S_G A_Y K^*.$$

\end{pr}
\begin{pre}
Consider the orthonormal basis $(\phi_n)_{n\geq1}$ for $\mathcal H(\Omega).$
For $v \in {\cal H}(\Omega),$ we have the following representation
$$
v=\sum_{k=1}^\infty(v,\phi_k)_{0,\Omega} \ \phi_k
$$
and since $\Gamma_0$ is bounded, we obtain that
\begin{eqnarray*}
\Gamma_0v&=&\sum_{k=1}^\infty(v,\phi_k)_{0,\Omega}\Gamma_0\phi_k\\
&=&\sum_{k=1}^\infty(v,\phi_k)_{0,\Omega}\kappa_ky_k\\
&=&\sum_{k=1}^\infty(v,\kappa_k\phi_k)_{0,\Omega}y_k \\
&=&\sum_{k=1}^\infty(v,\K y_k)_{0,\Omega}y_k\\
&=&S_Y\A_Y\K^*v,
\end{eqnarray*}
and that
\begin{eqnarray*}
\Gamma_0v&=&\sum_{k=1}^\infty(v,\phi_k)_{0,\Omega}\Gamma_0\phi_k\\
&=&\sum_{k=1}^\infty(v,\phi_k)_{0,\Omega}\kappa_ky_k\\
&=&\sum_{k=1}^\infty(v,\kappa_k\phi_k)_{0,\Omega}y_k\\
&=&\sum_{k=1}^\infty(v,\Gamma_0^*g_k)_{0,\Omega}y_k\\
&=&S_Y\A_G\Gamma_0v.
\end{eqnarray*}
Hence, $$\Ga_0= S_Y A_YK^*= S_Y A_G \Gamma_0.$$
In a similar way, one can establish that
$$K^* = S_G A_G \Gamma_0 = S_G A_Y K^*.\  \square $$ \\ 
\end{pre}
An interesting consequence of Proposition 4.1 is the following:
\begin{co}

The synthesis and analysis operators associated with the sequences  $(g_n)_{n\geq 1}$ and  $(y_n)_{n\geq 1}$ 
satisfie 
$$ S_Y A_G \subset I_{L^2(\p \Omega)} $$
and
$$ S_G A_Y \subset I_{L^2(\p \Omega)}.$$
\end{co}

\begin{pre}
From Proposition 4.1, we have $$\Ga_0= S_Y A_YK^*= S_Y A_G \Gamma_0.$$
Moreover, we showed in Lemma 4.2 that $\mathcal R(\Gamma_0)$ is dense in $L^2(\p \Omega).$ We therefore  have
$$ S_Y A_G \subset I_{L^2(\p \Omega)} .$$ Similarly, we have according to Proposition 4.1 that $$K^* = S_G A_G \Gamma_0 = S_G A_Y K^*,$$ and by  Corollary 4.8 that the operator $K$ is injective, which implies that $\mathcal R(K^*) $ is dense in $L^2(\p \Omega).$ Hence, $$ S_G A_Y \subset I_{L^2(\p \Omega)}. \ \square $$

\end{pre}
\begin{co}
The synthesis operators $S_G$ and $S_Y$ are bounded.
\end{co}
\begin{pre}
In Corollary 4.7, we showed that $(A_G^*,\mathcal D(A_G^*))$ and $(A_Y^*,\mathcal D(A_Y^*))$ are surjective, it follows that $\mathcal R(A_G)$ and $\mathcal R(A_Y)$ are closed. Moreover, from Corollary 4.9 the synthesis operators $(S_G,\mathcal D(S_G))$ and $(S_Y,\mathcal D(S_Y))$ are the inverse of the analysis operators $(A_Y,\mathcal D(A_Y))$ and  $(A_G,\mathcal D(A_G))$ respectively. Consequently, $S_G$ and $S_Y$ are bounded. $\square$
\end{pre} \\

Now, we can deduce the principal result of this paper.
\begin{co}
The sequences $(y_n)_{n\geq 1}$ and $(g_n)_{n\geq 1}$ are Riesz bases for $L^2(\p \Omega).$
\end{co}
\begin{pre}
 Having shown in Corollary 4.10 that the synthesis operators $S_G$ and $S_Y$ associated  with the sequences $(g_n)_{n\geq1}$ and $(y_n)_{n\geq1}$ respectively are bounded. It follows according to Lemma 2.4 that the sequences $(g_n)_{n\geq1}$ and $(y_n)_{n\geq1}$ are Bessel sequences. On the other hand, $(g_n)_{n\geq1}$ and $(y_n)_{n\geq1}$ are biorthogonal and complete according to Lemma 4.1 and Corollary 4.7. Consequently, $(g_n)_{n\geq1}$ and $(y_n)_{n\geq1}$ are Riesz Bases according to Theorem 2.2. $\square$
\end{pre}

\section{Regularity result for the Dirichlet problem}
Let $\Omega$ be a bounded Lipschitz domain in $\mathbb R^d.$ For $0\leq s \leq 1,$ we denote by $\mathcal H^s(\Omega)$ the space of real harmonic functions on the usual Sobolev space $H^s(\Omega)$, i.e., 
$$\mathcal H^s(\Omega) = \{ v\in H^s(\Omega)  \ | \ \Delta v=0 \ \mbox{in} \ \Omega \  \}$$ and  by $\mathcal H_s(\Omega)$ the following:
$$\mathcal H_s(\Omega) = \{ \  (I+F_1^*F_1)^{-s/2}v  \ | \ v \in \mathcal H(\Omega) \ \},$$
where $F_1$ is the Moore-Penrose inverse of the embedding operator $E_1.$
A first characterization of $\mathcal H^s(\Omega)$ is given in the following proposition:
 \begin{pr} Let $\Om \subset \mathbb R^d$ be a bounded Lipschitz domain.
 Then,  $\mathcal H^1(\Omega) = \mathcal H_1(\Omega)$ with an equivalence of norms.
   \end{pr}
\begin{pre} \
According to Proposition 2.4, we have $$E_1=(I+F_1^* F_1)^{-1/2}T_{F_1^*}$$
where $$T_{F_1^*}  = F_1^*(I+F_1F_1^*)^{-1/2} + E_1(I+F_1F_1^*)^{-1/2} .$$ 
Moreover, since the algebric equality $\mathcal{R}(E_1)=\mathcal H^1(\Omega)$ holds and that 
$$
\mathcal H_1(\Omega) = \{ \  (I+F_1^*F_1)^{-1/2} v  \ | \ v \in \mathcal H(\Omega) \ \},
$$

we obtain the algebric equality between $\mathcal H^1(\Om)$ and $\mathcal H_1(\Om),$ and all what is needed to prove is the equivalence of norms.
To this end,  for $u\in \mathcal H_1(\Om)$  consider the graph norm

$$\|u\|_{\mathcal H_1(\Om) }= \|(I+F_1^* F_1)^{1/2}   u\| _{0,\Om} .$$
For $v\in \mathcal H^1(\Om),$  $E_1 v \in \mathcal H_1(\Om)$ and
$$ \|(I+F_1^* F_1)^{1/2} E_1 v\| _{0,\Om} = \|T_{F_1^*} v\|_{0,\Om}. $$

Viewing $T_{F_1^*} $ is an isomorphism from $\mathcal H^1(\Om)$ into $\mathcal H(\Om)$ by Corollary 2.2. This assures the existence of  two positive  constants $c_1^{\prime}$ and $c_2^{\prime}$ not depending on $v$ such that
$$ c_1^{\prime}  \|v\|_{\p,\Om} \leq  \|(I+F_1^* F_1)^{1/2} E_1  v\| _{0,\Om}  \leq c_2^{\prime} \|v\|_{\p,\Om} .$$

Therefore, we deduce that the norms  $\|.\|_{\mathcal H_1(\Om)}$ and  $\|.\|_{\mathcal H^1(\Om)}$ are equivalent. $\square$
\end{pre}\\ \\
The following corollary could be deduced using the real interpolation method (see \cite{Ad}, \cite{Mcl} and \cite{T}).
   \begin{co} 
     Assume $0\leq s \leq 1,$ then $\mathcal H^s(\Omega)$ form an interpolatory family. Moreover, we have
     $$ \mathcal H^s(\Omega) =  \mathcal H_s(\Omega)$$ with equivalence of norms.
   \end{co}

 \begin{pr}
  Let $((\kappa_n, \phi_n))_{n\geq 1} $ be the sequence of couple in $\mathbb R_+ ^* \times \mathcal H(\Omega)$ associated with $\Gamma_0^* K^*$ such that  $\Gamma_0 K^* \phi_n = \kappa _n^2 \phi_n.$ Then the following statements are equivalent
  \begin{enumerate}
   \item $v\in \mathcal H^s(\Omega).$
   \item $\displaystyle{\sum _{n=1}^{+ \infty}}  \frac{1}{\kappa_n ^{4s}} \ |(v,\phi_n)_{0,\Omega}|^2 $ \mbox{converges.}
  \end{enumerate}
  \end{pr}
  \begin{pre}
  Let $v\in \mathcal H_s(\Omega).$  There exists then $\phi \in \mathcal H(\Omega)$ such that
  $$ v = (I+F_1^*F_1)^{-s/2}\phi ,$$
  which implies that
  $$(v,\phi_n)_{0,\Omega} = ((I+F_1^* F_1)^{-s/2}\phi , \phi_n)_{0,\Omega} = (\phi,(I+F_1^* F_1)^{-s/2}\phi_n)_{0,\Omega}.$$

 According to the Spectral Theorem  \cite [Theorem 5.1] {Co},  it follows that
   $$(v,\phi_n)_{0,\Omega} = ( \phi , (I+F_1^* F_1)^{-s/2}\phi_n)_{0,\Omega} = (\phi, \kappa_n^{2s} \phi_n),$$

which implies that
   $$ \frac{1}{\kappa_n^{2s}} (v,\phi_n)_{0,\Omega} = (\phi, \phi_n)_{0,\Omega}, $$

we therefore obtain that
   $$ \sum_{n=1}^{\infty} \frac{1}{\kappa_n ^ {4s}} |(v,\phi _n)_{0,\Omega}|^2 = \sum_{n=1}^{\infty} |(\phi,\phi_n)_{0,\Omega}|^2 . $$
   Since $\phi \in \mathcal H(\Om)$ and
   $ \displaystyle{\sum_{n=1}^{\infty}} |(\phi,\phi_n)_{0,\Omega}|^2 < + \infty, $ it follows then that

   $$  \displaystyle{\sum_{n=1}^{\infty}} \frac{1}{\kappa_n ^ {4s}} |(v,\phi _n)_{0,\Omega}|^2  < + \infty,  \ \ \  \forall v \in \mathcal H_s(\Omega),$$
   therefore according to Corollary 5.1, we conclude that  for $v\in \mathcal H^s(\Omega),$
    $  \displaystyle{\sum_{n}} \frac{1}{\kappa_n ^ {2s}} |(v,\phi _n)_{0,\Omega}|^2 $ converges.$\square$
   \end{pre}
   \\
   From Corollary 5.1 and Proposition 5.2, we deduce:
  
\begin{co}
We have $$ \mathcal R(K)\subset \mathcal R((I+F_1^*F_1)^{-1/4} P_{\mathcal H(\Om)})$$
and 
$$ \mathcal R(\Gamma_0^*)\subset \mathcal R((I+F_1^*F_1)^{-1/4} P_{\mathcal H(\Om)}).$$
\end{co}
\begin{pre}
Let $g \in L^2(\Omega).$ For all $n\geq 1,$ we have
$$ (Kg, \phi_n)_{0,\Omega}=(g,K^*\phi_n)_{0,\Omega}= (g, \kappa_n g_n)_{0, \p \Omega},$$
which leads to
$$ | \frac{1}{\kappa_n} (Kg, \phi_n)_{0,\Omega} | = |(g, g_n)_{0,\p\Omega} |$$
and since  $(g_n)_{n\geq 1}$ is a Bessel sequence, $\displaystyle{\sum_n} | (g,  g_n)_{0, \p \Omega} |^2$ converges,  and this  implies that $\displaystyle{\sum_n}| \frac{1}{\kappa_n} (Kg, \phi_n)_{0,\Omega} |^2$  converges, therefore $Kg\in \mathcal H ^{1/2}(\Omega).$ Similarly, for all $n\geq 1,$
$$ (\Gamma_0^*g, \phi_n)_{0,\Omega}=(g,\Gamma_0\phi_n)_{0,\Omega}= (g, \kappa_n y_n)_{0, \p \Omega}$$
so $$ | \frac{1}{\kappa_n} (\Gamma_0^*g, \phi_n)_{0,\Omega} | = |(g, y_n)_{0,\p\Omega} |.$$
Similarly, viewing  $(y_n)_{n\geq 1}$ is a Bessel sequence, it follows that
$\displaystyle{\sum_n} | (g,  y_n)_{0, \p \Omega} |^2$ converges and therefore $\displaystyle{\sum_n}| \frac{1}{\kappa_n} (\Gamma_0^*g, \phi_n)_{0,\Omega} |^2$  converges, thus $\Gamma_0^*g \in \mathcal H ^{1/2}(\Omega).$
\end{pre}
\\ \\
It turns out that the very weak solution of the Dirichlet problem for the Laplace equation (1), which is of the form $v= Kg$  lies in $\mathcal H^{1/2} (\Omega). \ \square $ \ \\

Now, we can state the main result of this section.
\begin{thm}
Let $\Omega \subset \mathbb R ^d$ be a bounded Lipschitz domain. For $g\in L^2(\partial\Omega)$, the very weak solution of the Dirichlet problem for the Laplace equation (1) lies in $H^\frac{1}{2}(\Omega)$ and there exist two positive constants $c_{\Omega}$ and $c_{\Omega}^{'}$  depending on  the geometry of $\Omega$ such that :
   $$
 c_{\Om}^{'} \ \|g\|_{L^2(\partial\Omega)}  \leq   \|v\|_{H^\frac{1}{2}(\Omega)}\leq c _{\Omega}\ \|g\|_{L^2(\partial\Omega)},
   $$ 
   Moreover, the solution operator $K$ is compact and injective.
\end{thm}

\begin{pre}
According to Corollary 5.2,  $$ \mathcal R(K)\subset \mathcal R((I+F_1^*F_1)^{-1/4} P_{\mathcal H(\Om)})$$
which implies that 
$$ \mathcal R(K)\subset \mathcal H^{1/2} (\Omega).$$ On the other hand, according to Lemma 2.7,  $\mathcal H^{1/2}(\Omega)$ is compactly imbedded in $L^2(\Omega)$, therefore $K$ is compact. The injectivity of $K$ holds from Corollary 4.8. On the other hand, since $$ \mathcal H_{1/2} (\Om)= \mathcal H^{1/2} (\Om),$$ and that
$$ \|Kg\|_{\mathcal H_{1/2} (\Om)} = \|(I+F_1^*F_1)^{1/4} Kg\|_{0,\Om},$$
we have
\begin{eqnarray*}
((I+F_1^*F_1)^{1/4} Kg,\phi_n)_{0,\Om} &=& (Kg, (I+F_1^*F_1)^{1/4} \phi_n)_{0,\Om}\\
&=& (Kg, \kappa_n^{-1} \phi_n)_{0,\Om} \\
&=& (g, \kappa_n^{-1} K^* \phi_n)_{0,\p \Om}\\
&=& (g,g_n)_{0,\p \Om}.
\end{eqnarray*}
Therefore, we obtain that $$ \|Kg\|_{\mathcal H_{1/2} (\Om)}^2 = \displaystyle{\sum_{n=1}^{\infty}} |(g,g_n)_{0,\p \Om} |^2 .$$
Moreover, since $(g_n)_{n\geq 1}$ is a Riesz basis for $L^2(\p \Om),$ there exist  two constants $b_G$ and $a_G$ such that 
$$
 a_G \|g\|_ {\p \Om}^2   \leq \|Kg\|_{\mathcal H_{1/2} (\Om)}^2  \leq  b_G \|g\|_{0,\p \Om}^2,$$

 where $b_G$ is the smallest upper bound of the sequence $(g_n)_{n\geq 1} $  and  $a_G$ is the greatest lower  bound of the sequence $(g_n)_{n\geq 1}, $ 
   therefore, we obtain that
$$\sqrt{a_G}\ \|g\| _{0,\p \Om} \leq   \|v\|_{\mathcal H_{1/2}(\Omega)}\leq \sqrt{b_G}\ \|g\|_{0,\partial\Omega}. \ \ 
     $$ 
Moreover, in view of Corollary 5.1, we have the equivalence of the norms $\|.\|_{\mathcal H_{1/2}(\Om)}$ and $\|.\|_{\mathcal H^{1/2}(\Om)},$ which implies that there exist two positive constants $c_{\Om}$ and $c_{\Om}^{\prime}$  depending on the geometry of $\Omega$ such that 
$$c_{\Om} ^{\prime} \ \|g\| _{0,\p \Om} \leq   \|v\|_{ H^{1/2}(\Omega)}\leq c_{\Om} \  \|g\|_{0,\partial\Omega}. \ \ 
    \square $$ 
\end{pre} \\ \\

{\bf{Remark: }}
 This work is part of the second author's ongoing  Ph.D. research, which is carried out at Moulay Ismail University, Meknes-Morocco.

\bibliographystyle{plain}

\end{document}